\documentclass[10pt]{article}
\usepackage{amsmath,latexsym,amsxtra,enumerate,hyperref,pdfsync,color,bbm,amsthm}
\usepackage[active]{srcltx}
\usepackage[english]{babel}
\theoremstyle{plain}
\allowdisplaybreaks
\newtheorem{lemma}{Lemma}[section]
\newtheorem{proposition}[lemma]{Proposition}
\newtheorem{theorem}[lemma]{Theorem}
\newtheorem{cor}[lemma]{Corollary}
\newtheorem{remark}[lemma]{Remark}
\newtheorem{definition}[lemma]{Definition}
\newtheorem{example}[lemma]{Example}
\newcommand{\sptext}[3]{\hspace{#1 em}\mbox{#2}\hspace{#3 em}}
\newcommand{\R}{\mathbbm{R}}
\newcommand{\Q}{\mathbbm{Q}}
\renewcommand{\P}{\mathbbm{P}}
\newcommand{\half}{\frac{1}{2}}
\newcommand{\ftn}{\mathcal{F}}
\newcommand{\gtn}{\mathcal{G}}
\newcommand{\vare}{\varepsilon}
\newcommand{\E}{\mathbbm{E}}
\newcommand{\pl}{\hspace{.1cm}}
\newcommand{\kla}{\left ( }
\newcommand{\mer}{\right ) }
\newcommand{\klae}{\left \{ }
\newcommand{\mere}{\right \} }
\newcommand{\equa}{\begin{eqnarray*}}
\newcommand{\tion}{\end{eqnarray*}}
\newcommand{\noo}{\left \|} 
\newcommand{\rrm}{\right \|} 
\newcommand{\bet}{\left |} 
\newcommand{\rag}{\right |} 
\newcommand{\ov}{\overline}
\newcommand{\nachop}[1]{\stackrel{#1}{\Longrightarrow}}
\begin{document}
\selectlanguage{english}

\title{Generalized fractional smoothness and $L_p$-variation of BSDEs
       with non-Lipschitz terminal condition}

\author{Christel Geiss, Stefan Geiss and 
        Emmanuel Gobet 
                       \\ \\
        Department of Mathematics  \\
        University of Innsbruck \\
        Technikerstra\ss e 13/7 \\
        A-6020 Innsbruck \\
        Austria \\
        christel.geiss@uibk.ac.at  and stefan.geiss@uibk.ac.at \\ \\
        Centre de Math\'ematiques Appliqu\'ees\footnote{A significant part of this work was done when the third author was affiliated to Grenoble Institute of Technology.}\\
        Ecole Polytechnique\\
        F-91128 Palaiseau Cedex\\
        France \\
        emmanuel.gobet@polytechnique.edu}
\maketitle

\begin{abstract}
We relate the $L_p$-variation, $2\le p < \infty$, of a solution of a backward stochastic 
differential equation with a path-dependent terminal condition to a generalized notion 
of fractional smoothness. This concept of fractional smoothness
takes into account the quantitative propagation of singularities in time.
\end{abstract}

MSC2010: 60H10, 46B70


\section*{Introduction}

During the last years the concept of fractional smoothness in the sense of function
spaces has been used in the theory of stochastic processes to analyze approximation
and variational properties. It turned out that phenomena known for special
examples can be explained in terms of fractional smoothness. For example, approximation
properties of certain stochastic integrals can be explained by the fractional smoothness
of the integral itself, see \cite{geiss:geiss:04,geiss-hujo}.
Similarly, variational properties of backward stochastic differential equations (BSDEs) can be upper
bounded in case that the fractional smoothness of the terminal condition is known.
To explain the latter aspect consider the BSDE
\[ Y_t = \xi + \int_t^T f(s,X_s,Y_s,Z_s) ds - \int_t^T Z_s dW_s \]
with a Lipschitz generator $f$, where $X=(X_t)_{t\in [0,T]}$ is a forward diffusion,
and define the $L_p$-variation
\[    {\rm var}_p (\xi,f,\tau) 
   :=   \sup_{i=1,...,n} \sup_{t_{i-1}< s \le t_i} \| Y_s - Y_{t_{i-1}} \|_p
      + \kla  \sum_{i=1}^{n} \int_{t_{i-1}}^{t_i} \| Z_t - \overline{Z}_{t_{i-1}} \|_p^2 dt\mer^\frac{1}{2} \]
where $\tau=(t_i)_{i=0}^n$ is a deterministic time-net $0=t_0<\cdots<t_n=T$,
\[   \overline{Z}_{t_{i-1}} 
   := \frac{1}{t_{i}-t_{i-1}} \E\left [ \int_{t_{i-1}}^{t_{i}} Z_s ds |\ftn_{t_{i-1}}
      \right ], \]
and where $2\le p < \infty$, which we will assume throughout this paper.
Note that by interchanging the $L_p$- and $L_2$-norms (where we use $p\ge 2$)
and using the Burkholder-Davis-Gundy inequality, the $L_p$-distance between the
stochastic integral $\int_0^T Z_s dW_s$ and its discrete counterpart 
$\sum_{i=1}^n \overline{Z}_{t_{i-1}} (W_{t_i}-W_{t_{i-1}})$ is upper bounded 
by a multiple of 
$\kla  \sum_{i=1}^{n} \int_{t_{i-1}}^{t_i} \| Z_t - \overline{Z}_{t_{i-1}} \|_p^2 dt\mer^\frac{1}{2}$.
Hence the quantity  ${\rm var}_p (\xi,f,\tau)$ is stronger compared to what is 
needed to quantify the discretization of the stochastic integral term of our BSDE.
Besides the fact that this variation gives a strong insight into the quantitative behavior
of the BSDE, in particular ${\rm var}_2 (\xi,f,\tau)$ was used to describe the error in adapted 
backward Euler schemes for $\xi = g(X_T)$ with $g$ being a Lipschitz function;
see 
\cite{Bouchard-Touzi,ZhangJ1} for implicit schemes and 
\cite{gobe:lemor:warin:2005,gobe:lemor:warin:2006} for explicit schemes possibly
with jump processes.
In \cite[Theorems 3.1 and 3.2]{gobe:makh:10} upper bounds for 
\[ \sum_{i=1}^{n} \int_{t_{i-1}}^{t_i} \| Z_t - \overline{Z}_{t_{i-1}} \|_2^2 dt \]
were obtained for $\xi = g(X_T)$ satisfying
\[ \E |g(X_T) - \E (g(X_T)|\ftn_t) |^2 \le c^2 (T-t)^\theta \]
for some $0<\theta \le 1$,
where $g$ is not assumed to be a Lipschitz function.
On the other hand, path-dependent settings without taking into account fractional
smoothness were considered, for example, in \cite{Hu-Ma:2004,Ma-Zhang:2002,ZhangY-ZhengW:2002}
In this paper, results are generalized and extended into the following
directions:

\begin{itemize}
\item We consider a path-dependent setting by terminal conditions of the form
      \[ \xi=g(X_{r_1},...,X_{r_L}) \]
      with $0=r_0<\cdots<r_L=T$, where $g$ is not necessarily a Lipschitz function
      and introduce a corresponding path-dependent fractional smoothness in the Malliavin sense.
      This concept of smoothness extends the classical concepts, based on real interpolation, to a time-dependent
      one taking care about the propagation of smoothness in time.
      In the classical case one would assign to a random variable $\xi$ some $0<\theta\le 1$
      which describes the fractional smoothness of $\xi$ while here we assign to the
      parameters $(\xi,f)$ of our BSDE a {\it vector} 
      $\Theta=(\theta_1,...,\theta_L)$, where $\theta_l$ stands for the local
      smoothness of the BSDE at time $r_l$. It turns out that this vector
      is completely characterized by the $L_p$-variation of $Y$ and $Z$.
      In case our terminal condition depends on $X_T$ only
      our generalized smoothness coincides with earlier approaches
      from, for example, \cite{geiss:geiss:04} and \cite{gobe:makh:10}.

\item Instead of the $L_2$-variation we consider the stronger $L_p$-variation with $2\le p < \infty$.
      In addition, the integrated $Z$-variation 
      $\sum_{i=1}^n \int_{t_{i-1}}^{t_i} \| Z_t - \overline {Z}_{t_{i-1}} \|_2^2 dt$ 
      is replaced by the variation $\| Z_s-Z_t\|_p$ with $s$ and $t$ being fixed, and the $L_p$-variation of the 
      process $Y$ is included as well.
      To our knowledge the weaker criterion for $0<p<2$ in the context of 
      this paper has not been considered yet and might require different arguments as some 
      of our proofs rely on the condition that $p\ge 2$.
 
\item We provide equivalences showing that the results are sharp. 

\item In Corollary \ref{cor:varition_with_theta_nets} we show, given the terminal condition 
      $\xi=g(X_{r_1},...,X_{r_L})$ has a certain fractional smoothness, how to obtain time-nets $\tau^n$ of
      cardinality $Ln+1$ such that
      \[  \sup_n \sqrt{n} \hspace{.2em} {\rm var}_p (\xi,f,\tau^n) < \infty. \]
      These time-nets compensate the possible singularities of the 
      $Z$-process when approaching a time-point $r_l$ from the left.
\end{itemize}

{\bf Organization of the paper.}
After introducing the setting in Section \ref{sec:setting},
we formulate in Section \ref{sec:fractional_smoothness_bsde} our concept of functional 
fractional smoothness of a BSDE and characterize this smoothness in various ways. Here we 
partly transfer the results from  \cite{geiss:geiss:04} and \cite{gobe:makh:10} from the
case $\xi=g(X_T)$ to the path-dependent one. 
In Section \ref{sec:sufficient_conditions}
we present two sufficient conditions for our fractional smoothness. The point of these 
two conditions (Corollary \ref{cor:sufficient_condition_I} and
Theorem \ref{thm:Lipschitz_functional}) is that they only involve the 
terminal condition $\xi$ and do not use the solution $Y$ nor the generator $f$ of our BSDE.
The proofs of the main results are contained in Section \ref{sec:proofs_of_main_results}.
\bigskip

{\bf Some notation.}
Given a vector $x\in\R^d$ we denote by $|x|$ its Euclidean norm,
for a linear operator $D\in {\cal L}(\R^n,\R^m)$ the symbol $|D|$ stands for 
the Hilbert-Schmidt norm, where $\R^n$ and $\R^m$ are equipped with the standard
Euclidean structure.
Given $D(t,x)\in {\cal L}(\R^n,\R^m)$ with
$(t,x)\in [0,T]\times \R^d$ and $0<T<\infty$, we use 
\[ \|D\|_\infty :=  \sup_{x\in\R^d,t\in [0,T]} |D(t,x)|. \]
Finally, 
$B(\eta_1,\eta_2) := \int_0^1 x^{\eta_1-1} (1-x)^{\eta_2-1} dx$
where 
$\eta_1,\eta_2>0$, 
will denote the Beta-function.


\section{Setting and basic concepts}
\label{sec:setting}

\subsection{Forward-backward stochastic differential equations}
\label{sec:forward-backward}

We fix a complete probability space $(\Omega,\ftn,\P)$, 
$T>0$, $d\ge 1$ and a $d$-dimensional standard Brownian motion $W=(W_t)_{t\in [0,T]}$ 
with $W_0\equiv 0$. Furthermore, we assume that 
$(\ftn_t)_{t\in [0,T]}$ is the augmentation of the natural filtration of $W$.

\paragraph{The forward equation.}
Let 
\[ X_t = x_0 + \int_0^t b(s,X_s)ds + \int_0^t \sigma(s,X_s) dW_s \]
with $x_0\in\R^d$, where $b:[0,T]\times\R^d\to \R^d$ and 
$\sigma:[0,T]\times \R^d \to {\cal L}(\R^d,\R^d)$ satisfy the following conditions:
\smallskip

\begin{enumerate}
\item[{\bf ($A_{b,\sigma}$)}]
      We have $b,\sigma \in C^{0,2}_b([0,T]\times \R^d)$, where the derivatives up to order two are taken with
      respect to the space-variables and, for some $\gamma \in (0,1]$, are assumed to be 
      $\gamma$-H\"older continuous (w.r.t. the parabolic metric) on all compact subsets of $[0,T]\times \R^d$.
      Moreover, there is a $\delta>0$ such that $\langle Ax,x \rangle \ge \delta |x|^2$ 
      for $x\in \R^d$ and
      $b$ and $\sigma$ are $\frac{ 1 }{2 }$-H\"older continuous in time, 
      uniformly in space.
\end{enumerate}
We work with the usual stochastic flow $(X_s^{t,x})_{s,t\in [0,T],x\in \R^d}$ 
that solves for $(t,x)\in [0,T]\times \R^d$ the SDE $X_s=x$ on $[0,t]$ and 
$dX_s^{t,x} = \sigma(s,X_s^{t,x})dW_s^t + b(s,X_s^{t,x}) ds$ on $[t,T]$, where
$W_s^t:= W_s-W_t$ and the augmented natural filtration $(\ftn_s^t)_{s\in [t,T]}$ 
of $(W_s^t)_{s\in [t,T]}$ is used
(i.p. $X=X^{0,x_0}$). With our assumptions
we can assume that $(X_s^{t,x})_{s,t\in [0,T], x\in \R^d}$ is a continuous process in $(s,t,x)$.
\medskip

If $g:\R^d\to \R$ is a polynomially bounded Borel function, $0<R\le T$, and
\begin{equation}\label{eqn:F}
   F(t,x) := \E g(X_R^{t,x})
   \sptext{1}{for}{1}
   0\le t \le R,
\end{equation}
then $F\in C^{1,2}([0,R)\times \R^d)$ and
\[   \frac{\partial}{\partial t} F(t,x)
   + \frac{1}{2} \left \langle A(t,x),
     D^2F(t,x) \right\rangle
   + \langle b(t,x),\nabla_x F(t,x)\rangle
   = 0 \]
by Proposition \ref{proposition:transition_density} below where
\[ D^2 :=  \left ( \frac{\partial^2}{\partial x_i \partial x_j} \right )_{i,j=1}^d. \]
The standard tail estimates for the transition density $\Gamma$  are re-called in 
Proposition \ref{proposition:transition_density}. They 
ensure that $\frac{ \partial }{ \partial_t} \nabla_x F$, 
$\nabla_x\frac{ \partial }{ \partial_t} F$ and $D^m_x F$ with $|m|\le 3$ exist and are continuous 
on $[0,R)\times \R^d$.
For $0\le t \le r < R\le T$ one has that, a.s.,
\equa
\nabla_x F(r,X_r^{t,x}) & = & \E \left ( g(X_R^{t,x}) N^{r,1,(t,x)}_R | \ftn_r^t \right ), \\
D^2  F(r,X_r^{t,x})     & = & \E \left ( g(X_R^{t,x}) N^{r,2,(t,x)}_R | \ftn_r^t \right )
\tion
for the Malliavin weights $N^{r,i,(t,x)}_R$ that satisfy, for any given $0<q<\infty$, that
\[ \left [ \E \left ( \left | N^{r,i,(t,x)}_R \right |^q | \ftn_r^t \right ) \right ]^\frac{1}{q}
   \le \frac{\kappa_q}{(R-r)^\frac{i}{2}} \hspace{.4em} a.s.
   \sptext{.7}{and}{.7}
   \E \left ( N^{r,i,(t,x)}_R | \ftn_r^t \right ) = 0 \hspace{.4em} a.s. \]
for $i=1,2$ and all $0\le t \le r < R \le T$ with a constant $\kappa_q>0$ independent from
$(t,r,R,x)$ 
(see \cite{gobe:munos:05}, \cite[Proof of Lemma 1.1]{gobe:makh:10} and
 Remark \ref{rem:Malliavin_weights} below). 
A typical application of these estimates are the crucial inequalities
\begin{eqnarray}
\|\nabla_x F(r,X_r^{t,x})\|_p &\le& \kappa_{p'} \frac{\|g(X_R^{t,x})-\E(g(X_R^{t,x}) | \ftn_r^t )\|_p}{\sqrt{R-r}}, \label{eqn:upper_bound_gradient}
\\
\|D^2 F(r,X_r^{t,x})\|_p & \le &   \kappa_{p'} \frac{\|g(X_R^{t,x})-\E(g(X_R^{t,x}) | \ftn_r^t )\|_p}{R-r}, \label{eqn:upper_bound_D2}
\end{eqnarray}
for $1<p,p'<\infty$ with $1=(1/p)+(1/p')$.
\medskip

\paragraph{The backward equation.}
We are interested in the backward equation
\[ Y_t  =  \xi + \int_t^T f(s,X_s,Y_s,Z_s) ds - \int_t^T Z_s dW_s
   \sptext{1.5}{for}{.5}
   t\in [0,T] \quad a.s. \]
and assume the following conditions:
\begin{enumerate}
\item[{\bf ($A_f$)}]
      The function $f:[0,T]\times \R^d \times \R \times \R^d \to \R$ is continuous in $(t,x,y,z)$ and
      continuously differentiable in $x$, $y$ and $z$ with uniformly bounded derivatives.
      In particular, 
      there are $K_f>0$ and $L_f>0$ such that
      \equa
       |f(s,x_1,y_1,z_1) - f(s,x_2,y_2,z_2)|\!\! &\le& \!\! L_f [|x_1-x_2| + |y_1-y_2|+|z_1-z_2|], \\
       |f(s,x,y,z)| \!\! &\le& \!\! K_f + L_f (|x|+|y|+|z|).
      \tion
\item[{\bf ($A_g$)}]
      There are $\mathcal{R}=\{r_0,...,r_L\}$ with
      $0=r_0<r_1<\cdots<r_L=T$ and a measurable function of at most polynomial growth
      $g:(\R^d)^L\to\R$ such that
      \[ \xi:= g(X_{r_1},...,X_{r_L}). \]
\end{enumerate}
In this setting, the solution $(Y,Z)$ to the above BSDE is uniquely defined in any $L_p$-space
for $1<p<\infty$; see \cite[Theorem 4.2]{Briand-et-al}. Additionally, we assume in the paper that the solution $(Y,Z)$ is 
realized such that, on $[r_{l-1},r_l)$,
\[ Y_t = u_l (\overline{X}_{l-1};t,X_t)
   \sptext{1}{and}{1}
   Z_t = v_l (\overline{X}_{l-1};t,X_t) \sigma(t,X_t),\]
 where we set $\overline{X}_{l-1}:=(X_{r_1}, \ldots, X_{r_{l-1}})$. The above functions $u_l$ and $v_l$ are well 
defined due to the next proposition, which is an extension of \cite[Theorem 3.2]{ZhangJ2} and
follows from Lemma  \ref{lemma:validity_of_assumptions}, see also \cite{Hu-Ma:2004}. 

\begin{proposition}\label{prop:Au}
Assume that $(A_{b,\sigma}), (A_f)$ and $(A_g)$ are satisfied.
Then, for $l=1,...,L$ there exist measurable $u_l: (\R^d)^{l-1}\times [r_{l-1},r_l)\times \R^d \to \R$
and
$v_l: (\R^d)^{l-1}\times [r_{l-1},r_l)\times \R^d \to \R^{1\times d}$ 
and Borel sets $D_l \subseteq \R^{d(l-1)}$, $l=2,...,L$, such that $D_l^c$ is of Lebesgue measure zero,
and such that
\begin{enumerate}[{\rm (i)}]
\item $u_l(\overline{x}_{l-1};\cdot,\cdot): [r_{l-1},r_l)\times \R^d \to \R$ is continuous and 
      continuously  differentiable w.r.t. the space variable with
      $\nabla_x u_l (\overline{x}_{l-1};t,x)=v_l (\overline{x}_{l-1};t,x)$, where
      $\overline{x}_{l-1}=(x_1, \ldots,x_{l-1})$,
\item there are $\alpha_l,q_{l,1},...,q_{l,l}\in [1,\infty)$ such that
      \begin{multline*}
            \sup_{t\in [r_{l-1},r_l)}|u_l(\overline{x}_{l-1};t,x)| 
          + \sup_{t\in [r_{l-1},r_l)} \sqrt{r_l-t} |v_l(\overline{x}_{l-1};t,x)| \\
      \le   \alpha_l(1+|x_1|^{q_{l,1}}+\cdots+|x_{l-1}|^{q_{l,l-1}}+|x|^{q_{l,l}}),
      \end{multline*}
\item for all $l=1,...,L$, $x_1,...,x_{l-1},x\in \R^d$ and $r_{l-1}\le s < r_l$ the triplet
      \[ \bigg (X_t^{s,x},u_l(\overline{x}_{l-1};t,X_t^{s,x}), 
                                      v_l(\overline{x}_{l-1};t,X_t^{s,x})
                                      \sigma(t,X_t^{s,x})
                                      \bigg )_{t\in [s,r_l)}
       \]
       solves the BSDE with generator $f$ and terminal condition
       \[ u_l(\overline{x}_{l-1};r_l,X_{r_l}^{s,x}) \] 
       where
       \[ u_l(\overline{x}_{l-1};r_l,x) :=
          \left \{ \begin{array}{rcl}
                   u_{l+1}(\overline{x}_{l-1},x;r_l,x) \chi_{D_l}(\overline{x}_{l-1})&:& 2\le l < L, \\
                   g(\overline{x}_{l-1},x) \chi_{D_l}(\overline{x}_{l-1})&:& l = L.
                   \end{array} \right .  \]
        and $u_1(r_1,x) := u_2(x;r_1,x)$.
\end{enumerate}
\end{proposition}
\bigskip
In the above proposition we used the convention that 
$h(\overline{x}_0;\cdot) := h(\cdot)$.
It should be noted that by Proposition \ref{prop:Au} 
we modify at each level $l=2,...,L$ the functional for the $Y$-process on a nullset. 
However, because of 
\begin{equation}\label{eqn:change_of_Y}
 \P(X_{r_1}\in D_2,...,(X_{r_1},...,X_{r_{L-1}})\in D_L)=1,
\end{equation}
this does not affect the $L_p$-solution of our BSDE so that Proposition \ref{prop:Au}
is sufficient for our purpose.

\paragraph{Piece-wise linearization of the backward equation.}
We let
 $F_l(\overline{x}_{l-1};\cdot,\cdot):[r_{l-1},r_l]\times \R^d \to \R$ be given by 
\[    F_l(x_1,...,x_{l-1};t,x) 
    = F_l(\ov{x}_{l-1};t,x)
   := \E u_l(x_1,...,x_{r_{l-1}};r_l,X_{r_l}^{t,x}). \]
The function $F_l$ solves the backward PDE
\[   \frac{\partial F_l}{\partial t}(\ov{x}_{l-1};t,x)
   + \frac{1}{2} \left \langle A(t,x),
     D^2F_l(\ov{x}_{l-1};t,x) \right\rangle
   + \langle b,\nabla_x F_l(\ov{x}_{l-1};t,x)\rangle
   = 0 \]
on the interval $[r_{l-1},r_l)$ for fixed $x_1,...,x_{l-1}\in\R^d$.

\paragraph{Two facts} that are frequently used in the paper.
Firstly, for a filtered probability space $(M,\Sigma,\Q,(\gtn_t)_{t\in [r,R]})$,
$1\le q \le \infty$, $r\le t \le R$ and $\xi\in L_q$, one has that
\begin{equation}\label{eqn:lp:cond:expect}
     \|\xi-\E(\xi|\gtn_t)\|_q
   \le
     \sup_{t \le s\le R}\|\xi-\E(\xi |\gtn_s)\|_q
   \le 2\|\xi-\E(\xi|\gtn_t)\|_q
\end{equation}
as a consequence that $\E(\cdot |\ftn_s)$ is a contraction on $L_q$.
Secondly, given the assumptions on our forward diffusion, a polynomially bounded
Borel function $g:\R^d\to \R$, $r\le t \le R \le T$ and $1\le q <\infty$, we have that
\begin{multline}\label{eqn:conditional_expectation_vs_transition_density}
    \noo g(X_R^{r,x}) - \E (g(X_R^{r,x})|\ftn_t^r) \rrm_q \\
\le \kla \int_{\R^d}\int_{\R^d}\int_{\R^d} |g(\xi)-g(\eta)|^q 
            \Gamma(r,x;t,y) \Gamma(t,y;R,\xi) \Gamma(t,y;R,\eta) d y d\xi d\eta 
            \mer^\frac{1}{q} \\
\le  2 \noo g(X_R^{r,x}) - \E (g(X_R^{r,x})|\ftn_t^r) \rrm_q.
\end{multline}


\subsection{Functional fractional smoothness}
\label{sec:fractional_smoothness_bsde}

The usage of fractional smoothness in the investigation of variational properties
of BSDEs is the central idea of this paper. Fractional smoothness can be defined in various
ways. One way is the so-called $K$-method, a  method where functions are decomposed into 
differentiable parts and parts that are not differentiable. A quantitative analysis of these
decompositions leads to fractional smoothness.
\smallskip

To be more precise, assume two Banach spaces $X_0$ and $X_1$, where (say) $X_1$ 
is continuously embedded into $X_0$, $0<t<\infty$ and $x\in X_0$, and recall that the 
$K$-functional is given by
\[     K(x,t;X_0,X_1) 
   := \inf \{ \|x_0\|_{X_0} + t \| x_1\|_{X_1}: x=x_0+x_1 \}.\]
For $0<\theta<1$ and $1\le q \le \infty$ this leads to the real interpolation spaces
\[    \| x \|_{(X_0,X_1)_{\theta,q}}
    := \noo t^{-\theta} K(x,t;X_0,X_1) \rrm_{L_q\left ((0,\infty),\frac{dt}{t}\right )} \]
with 
\[ X_1 \subseteq (X_0,X_1)_{\theta_1,q_1'}
       \subseteq (X_0,X_1)_{\theta_1,q_1}
       \subseteq (X_0,X_1)_{\theta_0,q_0}
       \subseteq X_0 \]
where $0<\theta_0 < \theta_1 < 1$ and $1 \le q_0,q_1,q_1' \le \infty$ with
$q_1' \le q_1$ (see \cite{Bennet-Sharpley,Bergh-Lofstrom}).
Applying this concept to the Malliavin Sobolev space $D_{1,p}$, we obtain the Malliavin Besov 
(or fractional Sobolev) spaces
\begin{equation}\label{eqn:besov-space}
 B_{p,q}^\theta := (L_p,D_{1,p})_{\theta,q}
\end{equation}
where $0<\theta<1$ is the main parameter of the smoothness and $1<q\le \infty$
the fine-tuning parameter. In a context close to this paper these spaces and related ones have been exploited for example in
\cite{geiss:geiss:04,geiss-hujo,gobe:makh:10,toivola:phd}. The classical setting of the
Wiener space is changed in \cite{geiss:geiss:04,gobe:makh:10} into a setting where the
standard Gaussian measure is replaced by the distribution of the forward diffusion.
Here we go one step ahead and replace $0<\theta<1$ by a vector $\Theta=(\theta_1,....,\theta_L)$, 
where $\theta_l$ describes the smoothness at time $r_l$:
\medskip
\begin{definition}\label{definition:functional_fractional_smoothness} \rm
Let $\Theta=(\theta_1,...,\theta_L)\in (0,1]^L$, $2\le p < \infty$
and $\xi\in L_p$. If $Y$ is the solution of the  
BSDE with generator $f$ and terminal condition $\xi$, then we let $(\xi,f)\in B_{p,\infty}^\Theta(X)$ provided that
there is some $c>0$ such that
\[ \| Y_{r_l}- \E(Y_{r_l}|\ftn_s) \|_p \le c (r_l-s)^\frac{\theta_l}{2}\]
for all $l=1,...,L$ and $r_{l-1}\le s < r_l$. The infimum over all possible $c>0$ is denoted
by 
\[ c_{B_{p,\infty}^\Theta}= c_{B_{p,\infty}^\Theta}(\xi,f).\]
In the case that $f=0$ we will simply write $\xi\in B_{p,\infty}^\Theta(X)$.
\end{definition}

\bigskip
Specializing to $p=2$ and to the linear one-step Gaussian case ($X=W$, $T=L=1$ and $f=0$)
it holds (see \cite[Corollary 2.3]{geiss-hujo}) that 
\[ g(W_1)\in B_{2,\infty}^{(\theta)}(W)
   \sptext{1}{if and only if}{1}
   g\in B_{2,\infty}^\theta(\R^d,\gamma_d),\]
where the Wiener space over the standard Gaussian measure $\gamma_d$ on $\R^d$ is considered.
In particular, for $d=1$ and for the orthonormal basis consisting of 
Hermite polynomials $(h_k)_{k=0}^\infty \subseteq L_2(\R,\gamma_1)$ we obtain that $g=\sum_{k=0}^\infty \alpha_k h_k\in B_{2,\infty}^\theta(\R,\gamma_1)$ if and only if there is some $c>0$ such that for all $0\le t< 1$ one has that
\[ \sum_{k=1}^\infty k t^{k-1} \alpha_k^2 \le \frac{c^2}{(1-t)^{1-\theta}}, \]
see \cite[Theorem 2.2]{geiss-hujo}.
These connections explain the notation $(p,\infty)$ in Definition 
\ref{definition:functional_fractional_smoothness}. For a more general connection
between the speed of convergence of the conditional expectations used in 
Definition \ref{definition:functional_fractional_smoothness} and the real interpolation method 
the reader is referred to \cite{geiss-hujo}.
Our definition of fractional smoothness by an upper bound of 
\[ \| Y_{r_l}- \E(Y_{r_l}|\ftn_s) \|_p \]
has the advantage that  
(\ref{eqn:upper_bound_gradient}) and (\ref{eqn:upper_bound_D2})
give the upper bounds
\[ \left \| \nabla_x F_l(\overline{X}_{l-1};s,X_s) \right \|_p
   \le \kappa_{p'} \frac{\| Y_{r_l}- \E(Y_{r_l}|\ftn_s) \|_p}{\sqrt{r_l-s}}
   \le \kappa_{p'} c_{B_{p,\infty}^\Theta} (r_l-s)^{\frac{\theta_l-1}{2}} \]
and
\[   \left \| D^2 F_l(\overline{X}_{l-1};s,X_s) \right \|_p
   \le \kappa_{p'} \frac{\| Y_{r_l}- \E(Y_{r_l}|\ftn_s) \|_p}{r_l-s} 
   \le \kappa_{p'} c_{B_{p,\infty}^\Theta} (r_l-s)^{\frac{\theta_l-2}{2}} \]
for $r_{l-1}\le s<r_l$ and $1=(1/p)+(1/p')$,
so that we can control the gradient and the Hessian of $F_l$.
\smallskip
For our paper the fine-tuning parameter $q=\infty$ in (the generalization of) (\ref{eqn:besov-space})
turns out to be the right one. 
\bigskip

Finally, we want to mention the coincidence, that most of the relevant examples are naturally linked to this fine-tuning parameter $q=\infty$ in (\ref{eqn:besov-space}).

\subsection{Time-nets, splines and entropy numbers}

In our BSDE system the $Z$-process gets possibly singular at any of the 
particular time points $r_l$ when $r_l$ is approached from the left. 
The degree of  this singularity is determined
by the parameter $\theta_l$ describing the fractional smoothness in $r_l$. To keep the variation
${\rm var}_p(g(X_{r_1},...,X_{r_L}),f,\tau)$ small, we have to choose time-nets which
refine on the left of $r_l$ with an order given by the fractional smoothness $\theta_l$ 
while each of the intervals $[r_{l-1},r_l]$ is divided into $n$ sub-intervals.

\medskip
\begin{definition}\rm
For $\Theta\in (0,1]^L$ we let $\tau^{n,\Theta}=(t_k^{n,\Theta})_{k=0}^{nL}$ be given by 
$ t_0^{n,\Theta}:=0$ and
\[    t_k^{n,\Theta}
   := r_{l-1} + (r_l-r_{l-1}) \left ( 1 - \left ( 1 - \frac{k-(l-1)n}{n} \right )^{\frac{1}{\theta_l}}
      \right )
   \sptext{.7}{for}{.7}
   (l-1)n < k \le ln.\]
\end{definition}
\bigskip

Estimates on  the $L_p$-variation $\|Y_t-Y_s\|_p$ are close
to estimates how good the process $Y$ can be approximated in $L_p$ by linear adapted splines, i.e. 
we simply compute adapted approximations
of $Y$ at the time-points $t_0,...,t_n$ and interpolate them linearly.
So the notion {\em adapted spline} refers to the fact that the knots are adapted, however the spline 
itself is not an adapted process.
The adapted splines are typically used in complexity theory for stochastic processes
to find efficient approximation schemes for stochastic processes where the whole 
path needs to be approximated but the adaptedness of the approximation is not fully needed,
see \cite{Creutzig-et-al}. Here we use the following notation:
\medskip

\begin{definition}\rm
Given a time-net $\tau=(t_k)_{k=0}^n$ with $r=t_0<\cdots<t_n=R\le T$ we say that the process $S=(S_t)_{t\in [r,R]}$ 
is an adapted spline based on $\tau$ provided that $S_{t_k}$ is $\ftn_{t_k}$-measurable for all
$k=0,...,n$ and
\[ S_t := \frac{t_k-t}{t_k-t_{k-1}} S_{t_{k-1}} + \frac{t-t_{t_{k-1}}}{t_k-t_{k-1}} S_{t_k}
   \sptext{1}{for}{1}
   t_{k-1}\le t \le t_k. \]
\end{definition}
\bigskip

Finally, we recall the notion of entropy numbers to measure and compare 
compactness properties of  $Y=(Y_s)_{s\in [t,r_l]}$ as $t\uparrow r_l$
where the process gets singular.

\begin{definition}\rm
Given a normed space $E$ and 
$A\subseteq E$ we define
$e_n(A|E) := \inf \vare$, where the infimum is taken over
all $\vare>0$ such that there are $x_1,...,x_n\in E$ with
\[ A \subseteq \bigcup_{i=1}^n \{ x_i + \vare B_E \} 
   \sptext{1}{with}{1}
   B_E : = \{ x\in E : \|x\| \le 1\}. \]
\end{definition}


\pagebreak
\section{Functional fractional smoothness and BSDEs}
\subsection{A general equivalence}

The basic result of this paper is  

\begin{theorem}
\label{thm:generalized_besov_bsde_local}
Assume that $(A_{b,\sigma})$, $(A_f)$ and $(A_g)$ are satisfied.
For $2\le p < \infty$ and {\em fixed} $l\in \klae 1,...,L\mere$ and $\theta_l \in (0,1]$ consider the following conditions:
\begin{enumerate}[{\rm ($C1_l$)}]
\item There is some  $c_1>0$ such that, for $r_{l-1}\le s<t<r_l$,
       \[ \| Z_t - Z_s \|_p \le c_1 \left ( \int_s^t (r_l-r)^{\theta_l-2} dr 
         \right )^\frac{1}{2}. \]
\item There is some $c_2>0$ with
      $\|Z_t\|_p \le c_2 (r_l-t)^{\frac{\theta_l-1}{2}}$
      for $r_{l-1} \le t < r_l$.
\item There is some $c_3>0$ such that, for $r_{l-1}\le s<t \le r_l$, 
      \[ \| Y_t-Y_s \|_p \le c_3 \left ( \int_s^t (r_l-r)^{\theta_l-1} dr 
         \right )^\frac{1}{2}.\]
\item There is some $c_4>0$ such that, for $r_{l-1}\le s< r_l$, 
      \[ \| Y_{r_l} - \E (Y_{r_l}|\ftn_s) \|_p \le c_4 (r_l-s)^\frac{\theta_l}{2}. \]
\item There is some $c_5>0$ such that, for $r_{l-1}\le t < r_l$,
      \[     \left \| \left ( 
             \int_{r_{l-1}}^t |(D^2 F_l)(\overline{X}_{l-1};s,X_s) 
                   |^2 ds 
             \right )^\frac{1}{2} \right \|_p
         \le c_5   (r_l-t)^\frac{\theta_l-1}{2} . \]
\item There is some $c_6>0$ such that for all $n=1,2,...$ there is an adapted spline $S^n=(S_t^n)_{t\in [r_{l-1},r_l]}$ 
      based on 
      \[ \kla r_{l-1} + (r_l-r_{l-1}) \left ( 1 - \left ( 1 - \frac{k}{n} \right )^{\frac{1}{\theta_l}}
              \right ) \mer_{k=0}^n \]
such that
      \[     \sqrt{n} \sup_{t\in [r_{l-1},r_l]} \| Y_t - S_t^n \|_p
         \le c_6. \]
      The spline can be arranged such that $S_{r_{l-1}}^n = Y_{r_{l-1}}$ and $S_{r_l}^n = Y_{r_l}$.
\item There is some $c_7>0$ such that for $r_{l-1}\le t < r_l$ one has that
      \[     \sup_{n\ge 1} \sqrt{n} e_n \big ( (Y_s)_{s\in [t,r_l]} | L_p \big )
         \le c_7 (r_l-t)^\frac{\theta_l}{2}.\]
\end{enumerate}
Then one has that  \quad
\begin{multline*}
   (C1_l)\!\nachop{\theta_l\in (0,1)} \!(C2_l) \Longleftrightarrow (C3_l) \Longleftrightarrow \\
   (C4_l) \Longleftrightarrow (C5_l) \Longleftrightarrow (C6_l) \Longleftrightarrow (C7_l)
        \Longrightarrow (C1_l). \end{multline*}
\end{theorem}

\begin{remark}\rm 
The implication  $(C1_l) \Longrightarrow (C2_l)$ does not hold in general. To see this we consider
$d=T=L=l=1$, $f=0$, $\theta_1=1$ and $p=2$, and let
\[ g = \sum_{n=0}^\infty \alpha_n h_n
   \sptext{1}{with}{1}
   \sum_{n=0}^\infty \alpha^2_n < \infty \]
where $(h_n)_{n=0}^\infty \subseteq L_2(\R,\gamma_1)$ is the orthonormal basis of 
Hermite polynomials. Then, as in \cite[Lemma 3.9]{geiss-hujo}, we get that
\[   \noo \frac{\partial^2 F_1}{\partial x^2}(t,W_t) \rrm_2^2 
   = \sum_{n=0}^\infty \alpha_{n+2}^2 (n+2)(n+1) t^n \] 
and
\[   \noo Z_t - Z_s \rrm_2^2 
   = \int_s^t \sum_{n=0}^\infty \alpha_{n+2}^2 (n+2)(n+1) r^n dr  .\]
Choosing $\alpha_n:=(n(n-1))^{-1/2}$ for $n\ge 2$ and $\alpha_0=\alpha_1=0$ gives 
($C1_l$) but $\sup_{0 \le t < 1} \| Z_t \|_2 = \infty$.
\end{remark}
\bigskip

From Theorem \ref{thm:generalized_besov_bsde_local} the multi-step case directly follows.
For its formulation we introduce 
for $\Theta=(\theta_1,...,\theta_L) \in (0,1]^L$ and $0\le t <T$ the function
\[    \varphi(t) 
   := \sum_{l=1}^L \chi_{[r_{l-1},r_l)} (t) 
      (r_l-t)^{\frac{\theta_l-1}{2}}. \]

\begin{theorem}
\label{thm:generalized_besov_bsde}
Assume that $(A_{b,\sigma})$, $(A_f)$ and $(A_g)$ are satisfied.
For $2\le p < \infty$ and $\Theta\in (0,1]^L$ consider the following conditions:
\begin{enumerate}
\item [{\rm (C1)}]
      There is some  $c_1>0$ such that, for $r_{l-1}\le s<t<r_l$,
       \[ \| Z_t - Z_s \|_p \le c_1 \left ( \int_s^t \frac{ \varphi(r)^2 }{ r_l-r} dr 
         \right )^\frac{1}{2}. \]
\item[{\rm (C2)}] There is some $c_2>0$ with
      $\|Z_t\|_p \le c_2 \varphi(t)$ for $0\le t <T$.
\item [{\rm (C3)}]There is some $c_3>0$ such that, for $r_{l-1}\le s<t \le r_l$, 
      \[ \| Y_t-Y_s \|_p \le c_3 \left ( \int_s^t \varphi(r)^2 dr 
         \right )^\frac{1}{2}.\]
\item [{\rm (C4)}]$(\xi,f)\in B_{p,\infty}^\Theta(X)$.
\item [{\rm (C6)}]There is some $c_6>0$ such that for all $n=1,2,...$ there is an adapted spline $S^n=(S_t^n)_{t\in [0,T]}$ 
      based on $\tau^{n,\Theta}$ such that
      \[     \sqrt{n} \sup_{t\in [0,T]} \| Y_t - S_t^n \|_p
         \le c_6. \]
\end{enumerate}
Then one has that  \quad
\[ (C1)\!\nachop{\Theta\in (0,1)^L} \!(C2) \Longleftrightarrow (C3) \Longleftrightarrow
   (C4) \Longleftrightarrow (C6) 
        \Longrightarrow (C1). \]
\end{theorem}
\bigskip

The remaining properties $(C5_l)$ and $C7_l)$ could be included as well.
By using the properties (C3) and (C1) we deduce by a simple computation

\begin{cor}\label{cor:varition_with_theta_nets}
For $0<\theta_l' < \theta_l < 1$, $l=1,...,L$ and $(\xi,f)\in  B_{p,\infty}^{\Theta}(X)$
one has that
\[ \sup_n \sqrt{n}\hspace{.2em} {\rm var}_p(\xi,f,\tau^{n,\Theta'}) < \infty. \] 
\end{cor}

Examples will be considered in 
Example \ref{example:bv} 
and Theorem \ref{thm:Lipschitz_functional}.
The proof of Theorem \ref{thm:generalized_besov_bsde_local} is postponed to 
Section \ref{subsection:proof:theo:besov}.
\medskip

\subsection{Sufficient conditions for fractional smoothness}
\label{sec:sufficient_conditions}

In this section we describe sufficient conditions on $\xi$ for the condition 
$(\xi,f)\in B_{p,\infty}^\Theta(X)$ which are {\em independent from the generator $f$}.
Note that in the case $L=1$ it follows by definition that 
$(\xi,0)\in B_{p,\infty}^\Theta(X)$ implies that $(\xi,f)\in B_{p,\infty}^\Theta(X)$.
To our knowledge it is open whether it still holds for $L>1$.

\subsubsection{The first sufficient condition}
\label{sec:first_sufficient_condition}
The first sufficient condition is based on the concept to measure the fractional smoothness
of a random variable on the Wiener space by mixing the underlying Gaussian structure 
with an independent copy and to look how sensitive the given random variable is
with respect to this operation (see, for example, \cite{Hirsch}).
In our setting this would correspond to comparing, for example,
$g(X_1)$ with $g(X_1^\eta)$ where $X_1^\eta$ is defined via a Brownian motion
$W_t^\eta:= \sqrt{1-\eta^2}W_t + \eta B_t$ with $B$ being a Brownian motion independent from $W$
and $0\le \eta \le 1$.
Because we have a time-dependent structure we extend this concept by allowing more general
operations with $W$ and its independent copy $B$.
\smallskip

Let us consider two independent $d$-dimensional Brownian motions $W$ and $B$ on the same 
complete probability space 
$(\Omega,\ftn,\P)$ starting in zero, and let us denote by $(\ftn^W_t)_{t\in [0,T]}$ 
(resp. $(\ftn^B_t)_{t\in [0,T]}$ and $(\ftn^{W,B}_t)_{t\in [0,T]}$) the $\P$-augmentation of the natural 
filtrations of $W$ (resp. $B$ and $(W,B)$).
For a measurable function $\eta:[0,T]\mapsto[-1,1]$ we define 
the standard $d$-dimensional $\ftn^{W,B}$-Brownian motion 
\[ W^\eta_t:=\int_0^t \sqrt{1-\eta(s)^2} dW_s+\int_0^t \eta(s) dB_s\]
and denote by $(\ftn^\eta_t)_{t\in [0,T]}$ the augmentation of its natural filtration.
We also define $X^\eta$ to be the strong $(\ftn_t^\eta)_{t\in [0,T]}$-measurable solution of
\[ X^\eta_t = x_0 + \int_0^t b(s,X^\eta_s)ds + \int_0^t \sigma(s,X^\eta_s) dW^\eta_s.\]
For a given $\ftn^\eta_T$-measurable terminal condition $\xi^\eta\in L_p$
with $2\le p <\infty$ we let $(Y^\eta,Z^\eta)$ be the $L_p$-solution in the filtration $(\ftn^\eta_t)_{t\in [0,T]}$ of 
\[  Y^\eta_t  =  \xi^\eta + \int_t^T f(s,X^\eta_s,Y^\eta_s,Z^\eta_s) ds - \int_t^T Z^\eta_s dW^\eta_s. \]
In the case $\eta\equiv 0$ we simply write $W=W^0$, $\xi=\xi^0$,
$(X,Y,Z)=(X^0,Y^0,Z^0)$, and $\ftn_t=\ftn_t^0$. Our aim is to bound the distance between 
$(X^\eta,Y^\eta,Z^\eta)$ and $(X,Y,Z)$ by  the following stability result:
\medskip

\begin{theorem}\label{theorem:sufficient_I_new}
Assume that $(A_{b,\sigma})$ and $(A_f)$ are satisfied.
Then for $2\le p<\infty$ and $\xi,\xi^\eta\in L_p$ we have that
\begin{multline*}
      \left \|\sup_{0\le t\le T}|X^\eta_t-X_t|\right \|_p
    + \left \|\sup_{0\le t\le T}|Y^\eta_t-Y_t|\right \|_p
    + \left \|\left (\int_0^T |Z^\eta_t-Z_t|^2 dt\right )^{1/2}\right \|_p \\
\le c \bigg [ \|\xi^\eta-\xi\|_p+ [1+\|\xi\|_p] \sqrt{\int_0^T \eta(t)^2 dt}\bigg ]
\end{multline*}
where $c>0$ depends at most  on $(p,T,b,\sigma,K_f,L_f)$ and 
is non-decreasing with respect to $K_f$ and $L_f$.
\end{theorem}
\bigskip

The proof can be found in Section \ref{subsection:proof:sec:sufficient_conditions}. 
The motivation for the result is Corollary \ref{cor:sufficient_condition_I} below.
To formulate it, given $0\le t < r \le T$ we let
\[ \eta_{t,r}(s) := \chi_{(t,r]}(s), \]
i.e. we replace the Brownian paths on $(t,r]$ by an independent copy.

\begin{cor}\label{cor:sufficient_condition_I}
Assume $2\le p <\infty$, ($A_{ b,\sigma }$), ($A_f$) and $\xi = g(X_{r_1},...,X_{r_L})\in L_p$ for some 
Borel measurable function $g:\R^L\to \R$. Let 
\[ \xi^{t,r}:=  g(X_{r_1}^{\eta_{t,r}},...,X_{r_L}^{\eta_{t,r}}) \]
for $0\le t < r \le T$ and  let $\Theta=(\theta_1,...,\theta_L)\in (0,1]^L$.
If there is a constant $c>0$ such that one has that
\begin{equation}\label{eqn:condition_Wiener_space}
      \| \xi - \xi^{t,r_l}\|_p
   \le c (r_l-t)^\frac{\theta_l}{2}
\end{equation}
for all $l=1,...,L$ and $r_{l-1} \le t<r_l$, then 
$(\xi,f) \in B_{p,\infty}^\Theta(X)$.
\end{cor}

\begin{proof}
For $r_{l-1} \le t<r_l$ we get by (\ref{eqn:conditional_expectation_vs_transition_density}) that 
\equa
      \| Y_{r_l} - \E(Y_{r_l}|\ftn_t) \|_p
&\le& \| Y_{r_l} - Y_{r_l}^{\eta_{t,r_l}} \|_p \\
&\le& c_{(\ref{theorem:sufficient_I_new})}
       \left [ 
            \| \xi - \xi^{t,r_l} \|_p 
          +  [1+\|\xi\|_p]\sqrt{\int_0^T \eta_{t,r_l}(r)^2 dr}
         \right ] \\
&\le& c_{(\ref{theorem:sufficient_I_new})}
       \left [ c  (r_l-t)^\frac{\theta_l}{2}
           +  [1+\|\xi\|_p] \sqrt{r_l-t}
         \right ].
\tion
\end{proof}

Using a truncation argument, we obtain a modified version of  Theorem \ref{thm:generalized_besov_bsde}, 
without assuming that $g$ is polynomially bounded nor that $f$ is continuously differentiable in $(x,y,z)$.

\begin{cor}\label{cor:sufficient_condition_I:new}
Assume ($A_{ b,\sigma }$) and that the generator  $f:[0,T]\times \R^d \times \R \times \R^d \to \R$ is 
continuous in $(t,x,y,z)$ and that  there is some $L_f>0$ such that
      \[     |f(s,x_1,y_1,z_1) - f(s,x_2,y_2,z_2)|
         \le L_f [|x_1-x_2| + |y_1-y_2|+|z_1-z_2|]. \]
Let $2\le p <\infty$, $\xi = g(X_{r_1},...,X_{r_L})\in L_p$ for some 
Borel measurable function $g:\R^L\to \R$, $\Theta\in (0,1]^L$ and let
$(Y,Z)$ be the $L_p$-solution of our BSDE.
Assume that condition (\ref{eqn:condition_Wiener_space}) is satisfied.
Then there are sets $\mathcal{N}_l \subseteq [r_{l-1},r_l)$ of Lebesgue measure zero such that the following is satisfied:

\begin{enumerate}
\item [{\rm (C1')}] There is some  $c_1>0$ such that for $s,t \in [r_{l-1},r_l) \setminus \mathcal{N}_l$ with  $r_{l-1}\le s<t<r_l$
      one has
      \[ \| Z_t - Z_s \|_p \le c_1 \left ( \int_s^t \frac{ \varphi(t)^2 }{ r_l-r} dr  
         \right )^\frac{1}{2}. \]
\item [{\rm (C2')}] There is some $c_2>0$ with
      $\|Z_t\|_p \le c_2 \varphi(t)$ for $t\in \bigcup_{l=1}^L ([r_{l-1},r_l) \setminus \mathcal{N}_l)$.
\item [{\rm (C3')}] There is some $c_3>0$ such that, for $r_{l-1}\le s<t \le r_l$, one has
      \[ \| Y_t-Y_s \|_p \le c_3 \left ( \int_s^t \varphi(r)^2 dr 
         \right )^\frac{1}{2}.\]
\end{enumerate}
\end{cor}
\medskip

\begin{proof}
(a) Let $(f^N)_{N\geq 1}$ be a sequence of generators satisfying assumption $(A_f)$ such that
\begin{enumerate}[(i)]
\item $\lim_N \left \|\int_0^T |f^N(s,X_s,Y_s,Z_s)-f(s,X_s,Y_s,Z_s)| ds\right \|_{p}=0$,
\item $K_{f^N}\leq 2 K_f$ and $L_{f^N}\leq L_f$. 
\end{enumerate}

(b) Letting $y^N=-N\lor y \land N$ for $y\in \R$ and $N\ge 1$,
$\xi^N$ satisfies ($A_g$) and $\|\xi^N-\xi\|_p\rightarrow 0$ as $N\rightarrow \infty$. 
In addition, for all $l=1,...,L$ and $r_{l-1} \le t<r_l$ we have
\[     \| \xi^N - (\xi^N)^{t,r_l}\|_p=\| \xi^N - (\xi^{t,r_l})^N\|_p\leq \| \xi - \xi^{t,r_l}\|_p
   \le c_{(\ref{eqn:condition_Wiener_space})} (r_l-t)^\frac{\theta_l}{2}.\]
(c) To $(\xi^N,f^N)$ we associate $(Y^N,Z^N)$ as BSDE solution in $L_p$.
In view of the inequality above and 
according to Corollary \ref{cor:sufficient_condition_I}, $(\xi^N,f^N) \in B_{p,\infty}^\Theta(X)$. 
Because $K_{f^N}$, $L_{ f^N }$ and $\|\xi^N \|_p$ are bounded independently of $N$, we 
have \[ \sup_{N\geq 1}c_{B_{p,\infty}^\Theta}(\xi^N,f^N)<\infty, \] 
which follows by the proof of Corollary \ref{cor:sufficient_condition_I}.
Theorem \ref{thm:generalized_besov_bsde} applies to $(Y^N,Z^N)$ for each $N$ and there are 
$c^N>0$ such that
\[ \| Z_{N,t} - Z_{N,s} \|_p \le c^N \left ( \int_s^t \frac{ \varphi(t)^2 }{ r_l-r} dr  \right )^\frac{1}{2} \]
for $r_{l-1} \le s < t < r_l$.
Looking at the constants in the proof of $(C4_l) \Rightarrow (C1_l)$ we realize that we can take
$\sup_N c^N =: c < \infty$.
By Lemma \ref{lemma:compare_briand_etal} applied to 
$\xi^{(0)}=\xi$, $f_0(\omega,s,y,z):= f(s,X_s(\omega),y,z)$, $(Y^{(0)},Z^{(0)})= (Y,Z)$, 
and 
$\xi^{(1)}=\xi^N$, $f_1(\omega,s,y,z):= f^N(s,X_s(\omega),y,z)$, 
$(Y^{(1)},Z^{(1)})= (Y^N,Z^N)$, 
there is a sub-sequence $(N_k)_{k=1}^\infty$ such that
$Z_{N_k,t}$ converges to $Z_t$ a.s. for $t\in [r_{l-1},r_l)\setminus \mathcal{N}_l$ for some  $\mathcal{N}_l$
of Lebesgue measure zero. Fatou's lemma gives
\[ \| Z_t - Z_s \|_p \le c \left ( \int_s^t \frac{ \varphi(t)^2 }{ r_l-r} dr  \right )^\frac{1}{2} \]
for $r_{l-1} \le s < t < r_l$ with $s,t\in [r_{l-1},r_l)\setminus \mathcal{N}_l$.
As in the proof of $(C1_l)\Longrightarrow (C2_l) \Longrightarrow (C3_l)$ below we can deduce $(C2')$ and $C3')$ where 
in the case $r_{l-1} \in \mathcal{N}_l$ in
$(C1_l)\Longrightarrow (C2_l)$ we have to replace $\| Z_{r_{l-1}}\|_p$ by  $\liminf_n \| Z_{\rho_n}\|_p$
with $\rho_n \in [r_{l-1},r_l)\setminus \mathcal{N}_l$ and $\rho_n \downarrow r_{l-1}$
\end{proof}

\begin{definition}
A measurable function $g:\R\to\R$ is of bounded variation,
in short $g\in {\rm BV}$, provided that 
\[  \| g \|_{BV}:= \sup_N \sup_{-\infty<x_0<\cdots < x_N<\infty} 
     \sum_{k=1}^N |g(x_k)-g(x_{k-1})| 
   < \infty. \]
\end{definition}

The following Example \ref{example:bv} is more general than needed in this paper, however 
this generality does not require any extra effort and constitutes the natural setting.

\pagebreak
\begin{example}\label{example:bv} 
Assume $0<\theta<\frac{1}{p}\le \alpha \le 1$,
$g_j\in {\rm BV}$ with $\sum_{j=1}^\infty \| g_j \|_{BV}^\alpha < \infty$, and
linear and continuous functionals $\mu_1,\mu_2, ... \in (C[0,T])^*$ 
with $\| \mu_j \| \le 1$ such that the laws of 
$\langle X,\mu_1 \rangle$, $\langle X,\mu_2 \rangle$, $\langle X,\mu_3 \rangle$, ... have densities bounded
uniformly by a constant $\beta >0$. Define 
\[ \xi := \Phi (g_1(\langle X, \mu_1 \rangle), g_2(\langle X, \mu_2 \rangle),... ), \]
where $\Phi$ is a measurable function such that
\[      | \Phi(x_1,x_2,...) - \Phi(y_1,y_2,...) | 
   \le \kappa \sum_{j=1}^\infty | x_j - y_j |^\alpha 
\]
for some $\kappa>0$. Then there is a constant $c>0$ such that for all measurable $\eta:[0,T]\to [-1,1]$ we have that
\[ \noo \xi - \xi^\eta \rrm_p \le c \left ( \int_0^T \eta(r)^2 dr \right )^\frac{\theta}{2}. \]
Consequently, given $\Theta \in (0,1/p)^L$ there is a constant $c' >0$ such that
\[ \| \xi - \xi^{t,r_l} \|_p \le c' (r_l-t)^{\frac{\theta_l}{2}} \]
for $r_{l-1} \le t < r_l$.
\end{example}
\smallskip

\begin{proof}
Using \cite[Theorem 2.4]{Avikainen_1} for $1 \le q < \infty$ we get that 
\equa
      \noo \xi - \xi^\eta \rrm_p
&\le& \kappa \noo \sum_{j=1}^\infty  |g_j(\langle X,\mu_j \rangle) - g_j(\langle X^\eta,\mu_j \rangle)|^\alpha
      \rrm_p \\
&\le& \kappa \sum_{j=1}^\infty  \noo g_j(\langle X,\mu_j \rangle) - g_j(\langle X^\eta,\mu_j \rangle)
                                     \rrm^\alpha_{\alpha p} \\
&\le& \kappa 3^{\alpha + \frac{1}{p}} \beta^{\frac{q}{q+1} \frac{1}{p}} 
      \sum_{j=1}^\infty \| g_j \|_{BV}^\alpha
                                         \noo 
                                         \langle X,\mu_j \rangle - \langle X^\eta,\mu_j \rangle
                                         \rrm_q^{\frac{q}{q+1} \frac{1}{p}} \\
&\le& \kappa 3^{\alpha + \frac{1}{p}} \beta^{\frac{q}{q+1} \frac{1}{p}} \sum_{j=1}^\infty \| g_j \|_{BV}^\alpha
      \sup_j \noo 
                                         \langle X,\mu_j \rangle - \langle X^\eta,\mu_j \rangle
                                         \rrm_q^{\frac{q}{q+1} \frac{1}{p}} \\
&\le& \kappa  3^{\alpha + \frac{1}{p}} \beta^{\frac{q}{q+1} \frac{1}{p}} \left [ \sum_{j=1}^\infty \| g_j \|_{BV}^\alpha \right ]
      \noo \sup_{0\le t \le T} |X_t-X^\eta_t| \rrm_q^{\frac{q}{q+1} \frac{1}{p}} \\
&\le& \kappa 3^{\alpha + \frac{1}{p}} \beta^{\frac{q}{q+1} \frac{1}{p}} \left [ \sum_{j=1}^\infty \| g_j \|_{BV}^\alpha \right ]
      \left ( c_{(\ref {eqn:estimate_difference_X})} \left ( \int_0^T \eta(r)^2 dr \right )^\frac{1}{2} \right )^{\frac{q}{q+1} \frac{1}{p}}\!\!\!\!,
\tion
where inequality $(\ref {eqn:estimate_difference_X})$ below is used.
Taking $1 \le q < \infty$ large enough the assertion follows.
\end{proof}

\subsubsection{The second sufficient condition}
The second sufficient condition relies on a simple 
iteration procedure: 
\medskip

\begin{theorem}\label{thm:Lipschitz_functional}
Assume that $(A_{b,\sigma})$ and $(A_f)$ are satisfied and that
\[ \xi := g(X_{r_1},...,X_{r_L}), \]
where 
\begin{multline*}
 \bet g(x_1,...,x_L) - g(x'_1,...,x'_L)\rag \\
\le \sum_{l=1}^L \left [ \bet g_l(x_l) - g_l(x'_l) \rag 
    + \psi_l(x_1,...,x_l; x'_1,...,x'_l) |x_l-x'_l| \right ]
\end{multline*}
with polynomially bounded Borel functions $g$, $g_l$ and $\psi_l$ such that
\begin{equation}\label{eqn:speed_of_convergence}
      \|g_l(X_{r_l})-\E(g_l(X_{r_l})|\ftn_t)\|_p 
   \le c (r_l-t)^\frac{\theta_l}{2}
\end{equation}
for $l=1,...,L$, $0<\theta_l \le 1$,  and $r_{l-1}\le t<r_l$.
Then, 
\[ (\xi,f) \in B_{p,\infty}^\Theta(X). \]
\end{theorem}

The proof of Theorem \ref{thm:Lipschitz_functional} is given in 
Section \ref{sec:proof_of_Theorem_thm:Lipschitz_functional}.

\bigskip
\begin{example}\rm
Let $\Phi:\R^L\to\R$ be Lipschitz and 
$g_1,...,g_L$ be as in Theorem \ref{thm:Lipschitz_functional}, and 
define
\[ 
  g(x_1,...,x_L) := \Phi(g_1(x_1),...,g_L(x_L)).
\]
\end{example}

To verify (\ref{eqn:speed_of_convergence}) for concrete functions $g_l$, it is sufficient
to check the inequality for the Brownian motion and for an appropriately rescaled
function:
\medskip

\begin{proposition}\label{proposition:change_gamma}
Let $c_{(\ref{proposition:transition_density})}>0$ be the constant from
Proposition \ref{proposition:transition_density} so that
\[            \Gamma (t,x;s,\xi)
          \le c_{(\ref{proposition:transition_density})} \pl  
              \gamma^{d}_{s-t}\kla\frac{x-\xi}{c_{(\ref{proposition:transition_density})}}\mer \]
and let $h_l(x):= g_l\kla x_0 + c_{(\ref{proposition:transition_density})} x \mer$ and
assume that
\begin{equation}\label{eqn:speed_of_convergence_gaussian}   
      \|h_l(W_{r_l})-\E(h_l(W_{r_l})|\ftn_t)\|_p 
   \le c_l (r_l-t)^\frac{\theta_l}{2}
   \sptext{1}{for}{1}
   0\le t < r_l,
\end{equation}
then (\ref{eqn:speed_of_convergence}) holds true for some $c>0$.
\end{proposition}

The proof of this proposition can be found in the appendix.
One can rescale the argument of the function $h_l$ in (\ref{eqn:speed_of_convergence_gaussian}) 
as well to assume that $r_l=1$.
Examples for (\ref{eqn:speed_of_convergence_gaussian}) with $d=1$ and $r_l=1$ are the following:

\begin{enumerate}[(a)]
\item If $h_l(x) = \chi_{[K,\infty)}(x)$ for some $K\in \R$, then
      $\theta = 1/p$ according to \cite[Example 4.7, Proposition 4.5]{geiss-toivola:1}.
\item If $h_l(x) = x^\alpha$ for $x  \ge 0$ and $h_l(x)=0$ otherwise, and $0<\alpha< 1-(1/p)$, 
      then $\theta= \alpha+(1/p)$ according to \cite[Example 5.2, Lemma 4.7]{toivola:1} and 
      \cite[Proposition 4.5]{geiss-toivola:1}.
\end{enumerate}

A precise investigation about the relation of (\ref{eqn:speed_of_convergence_gaussian}) 
to $B_{p,q}^\theta(\R^d,\gamma_d)$ can be found in \cite{geiss-toivola:2}.


\section{Proofs of the main results}
\label{sec:proofs_of_main_results}

\subsection{Proof of Theorem \ref{thm:generalized_besov_bsde_local}}
\label{subsection:proof:theo:besov}

\underline{$(C1_l) \Longrightarrow (C2_l)$} for $0<\theta_l<1$ is obvious as
\equa
      \| Z_t\|_p 
&\le&   \| Z_{r_{l-1}}\|_p 
      + c_1 \kla \int_{r_{l-1}}^t (r_l-r)^{\theta_l-2} dr \mer^\frac{1}{2} \\
& = & \| Z_{r_{l-1}}\|_p 
      + c_1 \kla \frac{1}{1-\theta_l}
            [(r_l-t)^{\theta_l-1} - (r_l-r_{l-1})^{\theta_l-1}]\mer^\frac{1}{2} \\
&\le& \| Z_{r_{l-1}}\|_p + c_1 (1-\theta_l)^{-\frac{1}{2}} (r_l-t)^\frac{\theta_l-1}{2}.
\tion
\underline{$(C2_l) \Longrightarrow (C3_l)$}  We observe that 
\equa
&   & \|Y_t - Y_s\|_p \\
& = & \left \|\int_s^t f(r,X_r,Y_r,Z_r) dr - \int_s^t Z_r dW_r\right \|_p \\
&\le& \int_s^t \|f(r,X_r,Y_r,Z_r)\|_p dr + 
      a_p \left (\int_s^t \|Z_r\|_p^2 dr\right )^\frac{1}{2} \\
&\le& K_f (t-s) +  L_f \int_s^t \| |X_r| + |Y_r| + |Z_r| \|_p dr
      + a_p \left (\int_s^t \|Z_r\|_p^2dr\right )^\frac{1}{2} \\
&\le& (t-s) \left [ K_f  + L_f \sup_{r\in [0,T]} \| X_r \|_p 
                         + L_f \sup_{r\in [0,T]} \| Y_r \|_p \right ] \\
&   & + c_2 (L_f \sqrt{T} + a_p) \left (\int_s^t (r_l-r)^{\theta_l-1} dr\right )^\frac{1}{2} 
\tion
where we used that $2\le p < \infty$ and where 
$a_p>0$ is the constant from the Burkholder-Davis-Gundy
inequality.
\bigskip\newline
\underline{$(C3_l) \Longrightarrow (C4_l)$} 
Here we get that 
\equa
       \| Y_{r_l} - \E( Y_{r_l}|\ftn_s) \|_p
&\le&  \| Y_{r_l} -     Y_s             \|_p + \| Y_s - \E( Y_{r_l}|\ftn_s) \|_p \\
&\le& 2 \| Y_{r_l} -     Y_s             \|_p \\
&\le& 2 c_3 \kla \int_s^{r_l} (r_l-r)^{\theta_l-1} dr \mer^\frac{1}{2} \\
& = & 2 c_3 \sqrt{\frac{1}{\theta_l}} (r_l-s)^\frac{\theta_l}{2}.
\tion
\underline{$(C4_l) \Longrightarrow (C5_l)$} 
We consider
\equa
&   & \left \| \left ( 
             \int_{r_{l-1}}^t |(D^2 F_l)(\overline{X}_{l-1};s,X_s) 
                   |^2 ds 
             \right )^\frac{1}{2} \right \|_p \\
& = & \left \| \left ( 
             \sum_{k=1}^d \int_{r_{l-1}}^t |(\nabla_x(\partial_{x_k} F_l))(\overline{X}_{l-1};s,X_s) 
                   |^2 ds 
             \right )^\frac{1}{2} \right \|_p \\
&\le& \frac{1}{\eta} \left \| \left ( 
             \sum_{k=1}^d \int_{r_{l-1}}^t |(\nabla_x(\partial_{x_k} F_l) \sigma)(\overline{X}_{l-1};s,X_s) 
                   |^2 ds 
             \right )^\frac{1}{2} \right \|_p \\
&\le& \sum_{k=1}^d \frac{1}{\eta} \left \| \left ( 
             \int_{r_{l-1}}^t |(\nabla_x(\partial_{x_k} F_l) \sigma)(\overline{X}_{l-1};s,X_s) 
                   |^2 ds 
             \right )^\frac{1}{2} \right \|_p \\
&\le&  \sum_{k=1}^d \frac{b_p}{\eta} \left \|
            \int_{r_{l-1}}^t (\nabla_x(\partial_{x_k} F_l) \sigma)(\overline{X}_{l-1};s,X_s) dW_s
                   \right \|_p
\tion
where $b_p>0$ is the constant from the Burkholder-Davis-Gundy inequality and 
the ellipticity condition on $\sigma$ implies that there exists an $\eta>0$ such that
\[ \eta |y|_{\R^d}  \le | y^* \sigma(t,x) |_{\R^d}
   \sptext{1}{for all}{1}
   x,y\in \R^d. \]
To upper-bound the terms of the last sum we use  It\^o's formula and our PDE (which reduces the 
number of terms) to obtain
\begin{eqnarray}
&   &   \partial_{x_k}F_l(\ov{X}_{l-1};t,X_t) 
      - \partial_{x_k}F_l (\ov{X}_{l-1};r_{l-1},X_{r_{l-1}}) \nonumber \\
& = & -\int_{r_{l-1}}^{t} \big [
           \langle \partial_{x_k}b,\nabla_x F_l \rangle 
         + \frac{1}{2}\langle \partial_{x_k}A,D^2 F_l \rangle
           \big ](\ov{X}_{l-1};s,X_s)ds \label{eqn:pde_derivative} \\
&   & +\int_{r_{l-1}}^{t}\left ( \nabla_x(\partial_{x_k}F_l)\sigma\right )
      (\ov{X}_{l-1};s,X_s)dW_s \nonumber
\end{eqnarray}
which implies that 
\equa
&   & \left \|
            \int_{r_{l-1}}^t (\nabla_x(\partial_{x_k} F_l) \sigma)(\overline{X}_{l-1};s,X_s) dW_s
                   \right \|_p \\
&\le& \|\nabla_x F_l(\ov{X}_{l-1};t,X_t)\|_p + \|\nabla_x F_l (\ov{X}_{l-1};r_{l-1},X_{r_{l-1}})\|_p\\
&   & + \noo \int_{r_{l-1}}^{t}\big [ \langle \partial_{x_k}b,\nabla_x F_l \rangle
      + \frac{1}{2} \langle \partial_{x_k}A,D^2F_l \rangle
        \big ](\ov{X}_{l-1};s,X_s)ds\rrm_p \\
&\le&   \kappa_{p'} \frac{R_t}{\sqrt{r_l-t}}
      + \kappa_{p'} \frac{R_{r_{l-1}}}{\sqrt{r_l-r_{l-1}}}
      + \kappa_{p'} \| \partial_{x_k} b\|_\infty \int_{r_{l-1}}^{r_l} \frac{R_s}{\sqrt{r_l-s}} ds\\
&   & + \kappa_{p'} \frac{\| \partial_{x_k} A\|_\infty}{2} \int_{r_{l-1}}^{r_l} 
        \frac{R_s}{r_l-s} ds
\tion
with
$R_s := \| Y_{r_l}- \E(Y_{r_l}|\ftn_s) \|_p$ and $r_{l-1}\le s < r_l$ where 
we used ($A_{b,\sigma}$) and inequalities (\ref{eqn:upper_bound_gradient}) and (\ref{eqn:upper_bound_D2}).
Consequently,
\equa
&   & \left \| \left ( 
             \int_{r_{l-1}}^t |(D^2 F_l)(\overline{X}_{l-1};s,X_s) 
                   |^2 ds 
             \right )^\frac{1}{2} \right \|_p \\
&\le& c_4 \frac{d b_p}{\eta} \kappa_{p'} \Bigg [   (r_l-t)^{\frac{\theta_l-1}{2}} 
                     + (r_l-r_{l-1})^{\frac{\theta_l-1}{2}} \\
&   & + \sup_{1\le k \le d} \| \partial_{x_k} b\|_\infty \int_{r_{l-1}}^{r_l} (r_l-s)^{\frac{\theta_l-1}{2}} ds \\
&   & + \sup_{1\le k \le d} \frac{\| \partial_{x_k} A\|_\infty}{2} \int_{r_{l-1}}^{r_l} 
                (r_l-s)^{\frac{\theta_l}{2}-1}  ds \Bigg ].
\tion
\underline{$(C5_l) \Longrightarrow (C2_l)$} 
Here we start with

\begin{lemma}\label{lemma:first-second-derivative}
Assume that $(A_{b,\sigma}), (A_f)$ and $(A_g)$ are satisfied.
There exists a constant $c>0$, depending at most on
$\sigma,b,T,d$ and $2 \le p <\infty$, such that, for all $r_{l-1} \le s < t < r_l$,
\equa
&   & \hspace*{-7em}
       \| \nabla_x F_l(\overline{X}_{l-1};t,X_t) - \nabla_x F_l(\overline{X}_{l-1};s,X_s) \|_p \\
&\le&   c (t-s) \|\nabla_x F_l(\overline{X}_{l-1};r_{l-1},X_{r_{l-1}}) \|_p \\
&   & + c (t-s) \noo \kla \int_{r_{l-1}}^s | D^2F_l (\overline{X}_{l-1};v,X_v)|^2 d v
                                                     \mer^\frac{1}{2}\rrm_p \\
&   & + c \noo \kla \int_s^t | D^2 F_l (\overline{X}_{l-1};v,X_v)|^2 d v
                    \mer^\frac{1}{2}\rrm_p.
\tion
\end{lemma}

\begin{proof}
For simplicity we will omit $\overline{X}_{l-1}$ in the computation.
Using (\ref{eqn:pde_derivative}) with $r_{l-1}$ replaced by $s$ 
we get that
\equa
&   & \hspace*{-3em}
      \|\nabla_x F_l(t,X_t) - \nabla_x F_l(s,X_s)  \|_p \\
&\le& \sum_{k=1}^d 
      \|\partial_{x_k} F_l(t,X_t) - \partial_{x_k} F_l(s,X_s)  \|_p \\
&\le& \left [ \sum_{k=1}^d \| \partial_{x_k} b \|_\infty \right ]
        \noo \int_s^t | \nabla_x F_l (v,X_v)| d v
                    \rrm_p \\
&   & + \left [ \sum_{k=1}^d \frac{\| \partial_{x_k} A\|_\infty}{2} \right ]
      \noo \int_s^t | D^2 F_l (v,X_v)| d v
                    \rrm_p \\
&   & + a_p \sum_{k=1}^d \noo \kla \int_s^t \bet (\nabla_x(\partial_{x_k}F_l)
      \sigma )(v,X_v)\rag^2 dv \mer^\frac{1}{2} \rrm_p \\
&\le& \left [ \sum_{k=1}^d \| \partial_{x_k} b \|_\infty \right ]
        \noo \int_s^t | \nabla_x F_l (v,X_v)| d v
                    \rrm_p \\
&   & + \left [ \sum_{k=1}^d \frac{\| \partial_{x_k} A\|_\infty (t-s)^\half}{2} 
                + d a_p  \| \sigma\|_\infty \right ]
      \noo \kla \int_s^t | D^2 F_l (v,X_v)|^2 d v \mer^\half
                    \rrm_p
\tion
where $a_p>0$ is the constant from the Burkholder-Davis-Gundy inequality, 
so that
\begin{multline}\label{eqn:preDelta_D2}
      \|\nabla_x F_l(t,X_t) -\nabla_x F_l(s,X_s) \|_p \\
\le    c_1 \noo \int_s^t | \nabla_x F_l (v,X_v)| d v
                   \rrm_p
     + c_2 \noo \kla \int_s^t | D^2 F_l (v,X_v)|^2 d v
                    \mer^\frac{1}{2}\rrm_p
\end{multline}
with
\[
c_1 :=  \sum_{k=1}^d \| \partial_{x_k} b \|_\infty
\sptext{1}{and}{1}
c_2 :=  \frac{\sqrt{T}}{2} \sum_{k=1}^d \| \partial_{x_k} A\|_\infty  + d a_p  \| \sigma\|_\infty. \]
Using this relation for $s=r_{l-1}$ and applying Gronwall's lemma
implies
\begin{multline*}
      \|\nabla_x F_l(t,X_t) \|_p \\
\le  e^{c_1 T} \left [       \|\nabla_x F_l(r_{l-1},X_{r_{l-1}}) \|_p
              + c_2 \noo \kla \int_{r_{l-1}}^t | D^2F_l (r,X_r) |^2 d r\mer^\frac{1}{2}
                    \rrm_p
      \right ].
\end{multline*}
Now we return to (\ref{eqn:preDelta_D2}) and get that 
\equa
&   &  \|\nabla_x F_l(t,X_t) - \nabla_x F_l(s,X_s) \|_p \\
&\le&   c_1 \int_s^t \| \nabla_x F_l(r,X_r) \|_p  dr
      + c_2 \noo \kla \int_s^t | D^2 F_l (r,X_r) |^2 d r\mer^\frac{1}{2}
            \rrm_p \\
&\le& \! c_1 e^{c_1 T} \! \int_s^t 
      \left [       \|\nabla_x F_l(r_{l-1},X_{r_{l-1}}) \|_p
              + c_2 \noo \kla \int_{r_{l-1}}^r | D^2 F_l (v,X_v)|^2 d v
                    \mer^\frac{1}{2}\rrm_p
      \right ] dr\\
&   & + c_2 \noo \kla \int_s^t | D^2 F_l (r,X_r) |^2 d r\mer^\frac{1}{2}
            \rrm_p \\
&\le&   c_1 e^{c_1 T} (t-s) \|\nabla_x F_l(r_{l-1},X_{r_{l-1}}) \|_p \\
&   & + c_1c_2 e^{c_1 T} \int_s^t 
              \noo \kla \int_{r_{l-1}}^s | D^2 F_l (v,X_v)|^2 d v
                    \mer^\frac{1}{2}\rrm_p dr \\
&   &      + c_1 c_2 e^{c_1 T} \int_s^t 
              \noo \kla \int_s^r | D^2 F_l (v,X_v)|^2 d v
                    \mer^\frac{1}{2}\rrm_p dr \\
&   & + c_2 \noo \kla \int_s^t | D^2 F_l (r,X_r) |^2 d r\mer^\frac{1}{2}
            \rrm_p  \\ 
&\le&   c_1 e^{c_1 T} (t-s) \|\nabla_x F_l (r_{l-1},X_{r_{l-1}}) \|_p \\
&   &  + (t-s)c_1 c_2 e^{c_1 T} 
              \noo \kla \int_{r_{l-1}}^s | D^2 F_l (v,X_v)|^2 d v
                    \mer^\frac{1}{2}\rrm_p \\
&   &      + [c_1 c_2 e^{c_1 T} (t-s) + c_2 ]
              \noo \kla \int_s^t | D^2 F_l (v,X_v)|^2 d v
                    \mer^\frac{1}{2}\rrm_p.
\tion
\end{proof}
For $r\in [r_{l-1},r_l)$ we consider
\begin{equation}\label{eqn:delta_vl}
     \delta v_l(\overline{x}_{l-1};r,x) 
  := v_l(\overline{x}_{l-1};r,x) - \nabla_x F_l (\overline{x}_{l-1};r,x)
\end{equation}
and get that, a.s.,
\begin{multline*}
u_l(\overline{x}_{l-1};r_{l-1},x_{l-1}) - F_l(\overline{x}_{l-1};r_{l-1},x_{l-1}) = \\
    \int_{r_{l-1}}^{r_l} \overline{f}(\overline{x}_{l-1};r,X_r^{r_{l-1},x_{l-1}}) dr \\
  - \int_{r_{l-1}}^{r_l} \delta v_l (\overline{x}_{l-1};r,X_r^{r_{l-1},x_{l-1}}) 
                         \sigma (r,X_r^{r_{l-1},x_{l-1}}) d W_r^{r_{l-1}}
\end{multline*}
with
\[    \overline{f}(\overline{x}_{l-1};r,x)
   := f(r,x,u_l(\overline{x}_{l-1};r,x),v_l(\overline{x}_{l-1};r,x)\sigma(r,x)). \]
Letting 
\[    \lambda^r(\overline{x}_{l-1};s,x) 
  :=  \int_{\R^d} \overline{f}(\overline{x}_{l-1};r,\xi) \nabla_x \Gamma(s,x;r,\xi) d\xi \]
and applying a stochastic Fubini argument, it follows that
\begin{multline*}
     \delta v_l(\overline{x}_{l-1};s,X_s^{r_{l-1},x_{l-1}}) \sigma(s,X_s^{r_{l-1},x_{l-1}}) \\
   = \int_s^{r_l} \lambda^r(\overline{x}_{l-1};s,X_s^{r_{l-1},x_{l-1}}) dr \pl
      {\sigma(s,X_s^{r_{l-1},x_{l-1}})}
      \quad a.s. 
\end{multline*}
for $s\in [r_{l-1},r_l)\setminus {\mathcal N}_l(\overline{x}_{l-1})$, where
${\mathcal N}_l(\overline{x}_{l-1})$ is a Borel set of measure zero. Hence for
$s\in [r_{l-1},r_l)\setminus {\mathcal N}_l(\overline{x}_{l-1})$ we get by
(\ref{eqn:upper_bound_gradient}) and Proposition \ref{prop:Au} that
\equa
&   & \| \delta v_l(\overline{x}_{l-1};s,X_s^{r_{l-1},x_{l-1}}) \sigma(s,X_s^{r_{l-1},x_{l-1}}) \|_p \\
&\le& \int_s^{r_l} \|\lambda^r(\ov{x}_{l-1};s,X_s^{r_{l-1},x_{l-1}})  \sigma(s,X_s^{r_{l-1},x_{l-1}})\|_p dr \\
&\le& \|\sigma \|_\infty \, \kappa_{p'} \int_s^{r_l} 
      \frac{\|f(\overline{x}_{l-1};r,X_r^{r_{l-1},x_{l-1}})\|_p}{\sqrt{r-s}} dr \\
&\le& \|\sigma \|_\infty  \kappa_{p'} 
      \int_s^{r_l} \bigg [ \frac{K_f + L_f \big [\|X_r^{r_{l-1},x_{l-1}}\|_p
      + \|u_l(\ov{x}_{l-1};r,X_r^{r_{l-1},x_{l-1}})\|_p \big ] }{\sqrt{r-s}} \\
&   &  \hspace*{10em}
       + \frac{ L_f 
       \|v_l(\ov{x}_{l-1};r,X_r^{r_{l-1},x_{l-1}})\sigma(r,X_r^{r_{l-1},x_{l-1}})\|_p}{\sqrt{r-s}} 
      \bigg ]
      dr \\
&\le& \|\sigma \|_\infty \kappa_{p'} \int_s^{r_l} 
      \frac{1} {\sqrt{r-s}}
      \bigg (
      K_f + L_f \bigg [\|X_r^{r_{l-1},x_{l-1}}\|_p + \\ 
&   & \hspace*{-1em}
        \alpha_l \left (  1 + \frac{\|\sigma \|_\infty }{\sqrt{r_l-r}} \right )
            \left \| 1+|x_1|^{q_{l,1}}+\cdots+|x_{l-1}|^{q_{l,l-1}} + |X_r^{r_{l-1},x_{l-1}}|^{q_{l,l}} \right \|_p \bigg ]
            \bigg )
           dr.
\tion
By continuity of both sides in $s$ one can estimate the first term by the last term in the 
above display for all $s\in [r_{l-1},r_l)$. Using the stochastic flow we obtain the inequality
\equa
&   & \| Z_s - \nabla_x F_l(\overline{X}_{l-1};s,X_s)\sigma(s,X_s) \|_p \\
&\le& \|\sigma \|_\infty \kappa_{p'} \int_s^{r_l} 
      \frac{1} {\sqrt{r-s}}
      \bigg (
      K_f + L_f \bigg [\|X_r\|_p + \\ 
&   &       \alpha_l \left (  1 + \frac{\|\sigma \|_\infty }{\sqrt{r_l-r}} \right )
            \left \| 1+|X_{r_1}|^{q_{l,1}}+\cdots+|X_{r_{l-1}}|^{q_{l,l-1}} + |X_r|^{q_{l,l}} \right \|_p \bigg ]
            \bigg )
           dr \\
&\le& c_0 <\infty
\tion
where $c_0>0$ does not depend on $s$.
The assertion ($C2_l$) follows from this and Lemma \ref{lemma:first-second-derivative} applied to $s=r_{l-1}$
because
\equa
      \| Z_r \|_p 
&\le& \| Z_r - \nabla_x F_l(\overline{X}_{l-1};r,X_r)\sigma(r,X_r) \|_p + \| \sigma \|_\infty \| \nabla_x F_l (\overline{X}_{l-1};r,X_r) \|_p \\
&\le& c_0 + \| \sigma \|_\infty (1+ c_{(\ref{lemma:first-second-derivative})}T) 
          \| \nabla_x F_l (\overline{X}_{l-1};r_{l-1},X_{r_{l-1}}) \|_p \\
&   &                  +\| \sigma \|_\infty c_{(\ref{lemma:first-second-derivative})}  \noo \kla \int_{r_{l-1}}^t | D^2 F_l (\overline{X}_{l-1};v,X_v)|^2 d v
                           \mer^\frac{1}{2}\rrm_p.
\tion

\underline{$(C4_l) \Longrightarrow (C1_l)$}
To make our assumption $(C4_l)$ more transparent, the constant $c_4>0$ of this condition is denoted
by $c_{B_{p,\infty}^\Theta}$ in the following. Using (\ref{eqn:delta_vl}) and letting
$r_{l-1}\le r < r_l$, by condition $(A_{b,\sigma})$ we get that 
\equa
&   &  \|  Z_t^{r_{l-1},\overline{x}_{l-1}} - Z_s^{r_{l-1},\overline{x}_{l-1}}  \|_p  \\
&\le&  \|  Z_t^{r_{l-1},\overline{x}_{l-1}} \sigma(t,X_t^{r_{l-1},x_{l-1}})^{-1} 
         - Z_s^{r_{l-1},\overline{x}_{l-1}} \sigma(s,X_s^{r_{l-1},x_{l-1}})^{-1}\|_p 
       \| \sigma\|_\infty\\
&   & + \| Z_s^{r_{l-1},\overline{x}_{l-1}}  
      \sigma(s,X_s^{r_{l-1},x_{l-1}})^{-1}(\sigma(t,X_t^{r_{l-1},x_{l-1}})- 
      \sigma(s,X_s^{r_{l-1},x_{l-1}}))\|_p \\
&\le& \|  \nabla_x F_l(\overline{x}_{l-1};t,X_t^{r_{l-1},x_{l-1}}) 
          - \nabla_x F_l(\overline{x}_{l-1};s,X_s^{r_{l-1},x_{l-1}}) \|_p\| \sigma\|_\infty \\
&   &  + \|    \delta v_l (\overline{x}_{l-1};t,X_t^{r_{l-1},x_{l-1}}) 
            -  \delta v_l (\overline{x}_{l-1};s,X_s^{r_{l-1},x_{l-1}}) 
            \|_p\| \sigma\|_\infty \\
&   & +  L_\sigma \|\sigma^{-1} \|_\infty \| Z_s^{r_{l-1},\overline{x}_{l-1}} \|_p \times \\
&   & (\| \E( \bet X_t^{r_{l-1},x_{l-1}} - X_s^{r_{l-1},x_{l-1}} \rag^p | \ftn_s^{r_{l-1}} )\|_\infty^\frac{1}{p} +|t-s|^\frac{1}{2})\\
&\le&    c_{\sigma,b,p,T} \left [ 
         D_1(\ov{x}_{l-1}) 
      +  D_2(\ov{x}_{l-1}) 
      +  D_3(\ov{x}_{l-1}) \right ]
\tion
with
\equa
D_1(\ov{x}_{l-1}) &:=&  \|  \nabla_x F_l(\overline{x}_{l-1};t,X_t^{r_{l-1},x_{l-1}}) 
          - \nabla_x F_l(\overline{x}_{l-1};s,X_s^{r_{l-1},x_{l-1}}) \|_p, \\
D_2(\ov{x}_{l-1}) &:=& \|  
               \delta v_l (\overline{x}_{l-1};t,X_t^{r_{l-1},x_{l-1}}) 
            -  \delta v_l (\overline{x}_{l-1};s,X_s^{r_{l-1},x_{l-1}})  \|_p, \\
D_3(\ov{x}_{l-1}) &:=&  (t-s)^\frac{1}{2} \| Z_s^{r_{l-1},\overline{x}_{l-1}} \|_p.
\tion
Now we show that each 
$\| D_i(\ov{X}_{l-1}) \|_p$ , $i=1,2,3$, is bounded by a constant times
$\kla \int_s^t (r_l-r)^{\theta_l-2} dr \mer^\frac{1}{2}$ which
implies $(C1_l)$.
\smallskip

\underline{The term $D_1(\ov{X}_{l-1})$:} Here we use Lemma \ref{lemma:first-second-derivative} 
to get
\equa
      \| D_1(\ov{X}_{l-1}) \|_p 
& = & \| \nabla_x F_l(\overline{X}_{l-1};t,X_t) - \nabla_x F_l(\overline{X}_{l-1};s,X_s) \|_p \\
&\le& c_{(\ref{lemma:first-second-derivative})}
      (t-s) \|\nabla_x F_l(\overline{X}_{l-1};r_{l-1},X_{r_{l-1}}) \|_p \\
&   & + c_{(\ref{lemma:first-second-derivative})}(t-s) \noo \kla \int_{r_{l-1}}^s | D^2F_l      
      (\overline{X}_{l-1};v,X_v)|^2 d v
                                                     \mer^\frac{1}{2}\rrm_p \\
&   & + c _{(\ref{lemma:first-second-derivative})}\noo \kla \int_s^t | D^2 F_l 
      (\overline{X}_{l-1};v,X_v)|^2 d v
                    \mer^\frac{1}{2}\rrm_p \\
&\le& c_{(\ref{lemma:first-second-derivative})}
      (t-s) \|\nabla_x F_l(\overline{X}_{l-1};r_{l-1},X_{r_{l-1}}) \|_p \\
&   & + c_{(\ref{lemma:first-second-derivative})}(t-s) 
      \kla \int_{r_{l-1}}^s  \| D^2F_l (\overline{X}_{l-1};v,X_v)\|_p^2 d v
                                                     \mer^\frac{1}{2} \\
&   & + c _{(\ref{lemma:first-second-derivative})}
      \kla \int_s^t \| D^2 F_l (\overline{X}_{l-1};v,X_v)\|_p^2 d v
                    \mer^\frac{1}{2} \\
&\le& c_{(\ref{lemma:first-second-derivative})} 
      (t-s)  \kappa_{p'} c_{B_{p,\infty}^\Theta} (r_l-r_{l-1})^{\frac{\theta_l-1}{2}}  \\
&   & + c_{(\ref{lemma:first-second-derivative})}(t-s) 
      \kla \int_{r_{l-1}}^s  
           \kappa_{p'}^2 c_{B_{p,\infty}^\Theta}^2 (r_l-v)^{\theta_l-2} d v
                                                     \mer^\frac{1}{2} \\
&   & + c _{(\ref{lemma:first-second-derivative})}
      \kla \int_s^t \kappa_{p'}^2c_{B_{p,\infty}^\Theta}^2 (r_l-v)^{\theta_l-2}d v
                    \mer^\frac{1}{2}
\tion
where we have used (\ref{eqn:upper_bound_D2}).
Finally we apply
\equa
       (t-s)      \kla \int_{r_{l-1}}^s   (r_l-v)^{\theta_l-2} d v \mer^\frac{1}{2}
 &\le& (t-s) \sqrt{s-r_{l-1}} (r_l-s)^{\frac{\theta_l-2}{2}} \\
 &\le& \sqrt{t-s} \sqrt{s-r_{l-1}} \kla \int_s^t   (r_l-v)^{\theta_l-2} d v \mer^\frac{1}{2}.
\tion
\underline{The term $D_2(\ov{x}_{l-1})$ and a linearization:}
First we follow the approach of \cite{gobe:makh:10} done for the one-step scheme,
that shows that the difference process $((v_l -\nabla_x F_l)(\overline{X}_{l-1};r,X_r))_{r\in [r_{l-1},r_l)}$
solves the  linear BSDE with the generator $f^{\rm lin}$ defined below. We fix $x_1,...,x_{l-1}\in\R^d$
and define
$f^{\rm lin}: [r_{l-1},r_l)\times \R^d \times \R^{1\times d} \times \R^{d\times d} 
  \to \R^{1\times d}$
by
\[ f^{\rm lin}(\overline{x}_{l-1};r,x,U,V):=A_l^0(\overline{x}_{l-1};r,x) 
      + U B_l^0(\overline{x}_{l-1};r,x) 
      + \sum_{j=1}^d V_j C_l^{j,0}(\overline{x}_{l-1};r,x), \]
where $V_j$ is the $j$-th row of $V$, with
\begin{multline*}
      A_l^0(\overline{x}_{l-1};r,x)
  :=  \nabla_x f \Big (r,x,u_l(\overline{x}_{l-1};r,x),
                          v_l(\overline{x}_{l-1};r,x) \sigma(r,x)\Big ) \\
    +  \frac{\partial f}{\partial y} \Big (r,x,u_l(\overline{x}_{l-1};r,x),v_l(\overline{x}_{l-1};r,x)\sigma(r,x)\Big )  
            \nabla_x F_l(\overline{x}_{l-1};r,x)\\
    + \sum_{j=1}^d \frac{\partial f}{\partial z_j} \Big (r,x,u_l(\overline{x}_{l-1};r,x),
                                                             v_l(\overline{x}_{l-1};r,x) \sigma(r,x)\Big ) \times \\
     \times  \nabla_x \kla \sum_{k=1}^d \frac{\partial F_l}{\partial x_k}(\overline{x}_{l-1};r,x) \sigma_{kj}(r,x) \mer,
\end{multline*}
\vspace*{-2em}
\begin{multline*}
      B_l^0(\overline{x}_{l-1};r,x) 
:=  \frac{\partial f}{\partial y}(r,x, u_l(\overline{x}_{l-1};r,x), v_l(\overline{x}_{l-1};r,x)\sigma(r,x)) 
      I_{\R^d}  \\
   + \nabla_x b(r,x)  
      + \sum_{j=1}^d \frac{\partial f}{\partial z_j}(r,x,u_l(\overline{x}_  
        {l-1};r,x),v_l(\overline{x}_{l-1};r,x)\sigma(r,x)) \nabla_x \sigma_j(r,x)
\end{multline*}
and
\begin{multline*}
      C_l^{j,0}(\overline{x}_{l-1};r,x) \\
:=  \frac{\partial f}{\partial z_j}(r,x,u_l(\overline{x}_{l-1};r,x),  v_l(\overline{x}_{l-1};r,x) 
      \sigma(r,x)) I_{\R^d} 
    + \nabla_x \sigma_j(r,x),
\end{multline*}
with $\sigma_j=(\sigma_{kj})_{k=1}^d\in \R^d$, $\delta v_l$ defined as in (\ref{eqn:delta_vl}), and 
\[ \delta u_l(\overline{x}_{l-1};r,x) := u_l(\overline{x}_{l-1};r,x)-F_l(\overline{x}_{l-1};r,x).\]
This implies
\begin{equation}\label{eqn:upper_bound_f_linear}
    |f^{\rm lin}(\overline{x}_{l-1};r,x,u,v)|
\le |A_l^0(\overline{x}_{l-1};r,x)| + c_{(\ref{eqn:upper_bound_f_linear})}
    [|u|+|v|].
\end{equation}
To associate a BSDE to the driver $f^{\rm lin}$, we first check that 
\begin{equation}
\label{label:flin}
\int_{r_{ l-1 }}^{r_l}\|A_l^0(\overline{x}_{l-1};r,X_r^{r_{l-1},x_{l-1}})\|_p \ dr <\infty.
\end{equation}
For this purpose we let
\[     \psi_l(\ov{x}_{l-1};r) 
   :=  1 + \| \nabla_x F_l (\overline{x}_{l-1};r,X_r^{r_{l-1},x_{l-1}}) \|_p
         + \|      D^2 F_l (\overline{x}_{l-1};r,X_r^{r_{l-1},x_{l-1}}) \|_p, \]
which implies that
\begin{equation}\label{eqn:upper_bound_Al0}
    \| A_l^0(\overline{x}_{l-1};r,X_r^{r_{l-1},x_{l-1}}) \|_p
\le c_{(\ref{eqn:upper_bound_Al0})} \psi_l(\ov{x}_{l-1};r).
\end{equation}
In view of 
(\ref{eqn:upper_bound_gradient}) and
(\ref{eqn:upper_bound_D2})
we have that
\begin{multline}\label{eqn:upper_bound_psi_l}
      \psi_l(\ov{x}_{l-1};r) 
 \le  1 + (1+\sqrt{r_l-r})\frac{\kappa_{p'}}{r_l-r}\times \\
      \hspace{7em} \left \| F_l(\ov{x}_{l-1};r_l,X_{r_l}^{r_{l-1},x_{l-1}}) -
                                F_l(\ov{x}_{l-1};r  ,X_r    ^{r_{l-1},x_{l-1}}) \right \|_p.
\end{multline}
To obtain the integrability of the upper bound on $\psi_l(\ov{x}_{l-1};r)$ (and thus that 
of $\|A_l^0(\overline{x}_{l-1};r,X_r^{r_{l-1},x_{l-1}})\|_p$), we show that the assumption 
on global fractional smoothness implies a local fractional smoothness. Indeed, our global 
assumption reads as
\begin{equation}\label{eqn:smoothness_of_F_l}
       \noo F_l(\ov{X}_{l-1};r_l,X_{r_l}) - F_l(\ov{X}_{l-1};s,X_s) \rrm_p
   \le c_{B_{p,\infty}^\Theta} (r_l-s)^\frac{\theta_l}{2}
\end{equation}
for $r_{l-1}\le s < r_l$.
For any $0<\delta<1$ this implies that
\[     \int_{r_{l-1}}^{r_l} (r_l-s)^{-\frac{p\theta_l}{2}-\delta}
       \noo F_l(\ov{X}_{l-1};r_l,X_{r_l}) - F_l(\ov{X}_{l-1};s,X_s) \rrm_p^p ds 
   < \infty. \]
Using the transition density of $X$ and Fubini's theorem implies the existence
of a Borel set $E_l\subseteq (\R^d)^{l-1}$ such that $E_l^c$ has Lebesgue measure zero
and
\[     \int_{r_{l-1}}^{r_l} (r_l-s)^{-\frac{p\theta_l}{2}-\delta}
       \noo F_l(\ov{x}_{l-1};r_l,X_{r_l}^{r_{l-1},x_{l-1}}) - F_l(\ov{x}_{l-1};s,X_s^{r_{l-1},x_{l-1}}) \rrm_p^p ds 
   < \infty \]
for all $(x_1,...,x_{l-1})\in E_l$. For those $(x_1,...,x_{l-1})\in E_l$ we may deduce 
(using (\ref{eqn:lp:cond:expect}))
for $s\in ( (r_{l-1}+r_l)/2,r_l)$ and $a_l:= s-(r_l-s)$ that 
\equa
&   & \noo F_l(\ov{x}_{l-1};r_l,X_{r_l}^{r_{l-1},x_{l-1}}) - F_l(\ov{x}_{l-1};s,X_s^{r_{l-1},x_{l-1}}) 
      \rrm_p^p \\
&\le& 2^p (s-a_l)^{-1}  (r_l-a_l)^{\delta + \frac{p\theta_l}{2}} \\
&   & \int_{a_l}^s (r_l-r)^{-\frac{p\theta_l}{2}-\delta} 
      \noo F_l(\ov{x}_{l-1};r_l,X_{r_l}^{r_{l-1},x_{l-1}}) - F_l(\ov{x}_{l-1};r,X_r^{r_{l-1},x_{l-1}}) \rrm_p^p dr \\
&\le& 2^{p + \delta + \frac{p\theta_l}{2}} (r_l-s)^{\delta + \frac{p\theta_l}{2}-1}  \\
&   & \int^{r_l}_{r_{l-1}} (r_l-r)^{-\frac{p\theta_l}{2}-\delta} 
      \noo F_l(\ov{x}_{l-1};r_l,X_{r_l}^{r_{l-1},x_{l-1}}) - F_l(\ov{x}_{l-1};r,X_r^{r_{l-1},x_{l-1}}) \rrm_p^p dr.
\tion
Taking $0<\delta<1$ such that $\delta + \frac{p\theta_l}{2}-1>0$ we obtain a local
fractional smoothness for all $(x_1,...,x_{l-1})\in E_l$.  Then for $\overline{x}_{l-1} \in E_l$
the inequality \eqref{label:flin} is satisfied.
Thus, because of \cite[Theorem 2.1]{gobe:makh:10} the process 
$(\delta v_l (\overline{x}_{l-1};s,X_s^{r_{l-1},x_{l-1}}))_{s\in [r_{l-1},r_l)}$
solves the $U$-component of the BSDE
\equa
      U_s^{r_{l-1},\ov{x}_{l-1}} 
& = & \int_s^{r_l} f^{\rm lin}(\overline{x}_{l-1};r,X_r^{r_{l-1},x_{l-1}},U_r^{r_{l-1},\ov{x}_{l-1}},
                     V_r^{r_{l-1},\ov{x}_{l-1}})dr \\
&   & - \kla \int_s^{r_l} (V_r^{r_{l-1},\ov{x}_{l-1}})^* dW_r^{r_{l-1}}\mer^*
\tion
for all $\overline{x}_{l-1}\in E_l$ (according to (\ref{eqn:upper_bound_f_linear}), (\ref{label:flin}) and
\cite[Theorem 4.2]{Briand-et-al} this BSDE has a unique $L_p$-solution).
\smallskip

\underline{Upper bound for $\|D_2(\ov{X}_{l-1})\|_p$:}
Applying Lemma \ref{lemma:upper_bound_V} to $h=f^{\rm lin}$ (the function $\kappa$
from Lemma \ref{lemma:upper_bound_V}(iii)
is obtained by Proposition \ref{proposition:transition_density} and (\ref{label:flin}) is used)
it follows that

\equa 
&   & \| U_s^{r_{l-1},\ov{x}_{l-1}} \|_p \\
&\le& c_{(\ref{lemma:upper_bound_V})}
       \left \| \int_s^{r_l}  |A_l^0(\ov{x}_{l-1};r,X_r^{r_{l-1},x_{l-1}})| dr \right \|_p \\
&\le& c_{(\ref{lemma:upper_bound_V})}  c_{(\ref{eqn:upper_bound_Al0})} 
      \int_s^{r_l} \psi_l(\ov{x}_{l-1};r)   dr \\
&\le& c_{(\ref{lemma:upper_bound_V})}  c_{(\ref{eqn:upper_bound_Al0})} 
      \bigg [  [r_l-s]+ \kappa_{p'} \int_s^{r_l} (1+\sqrt{r_l-r}) \\
&   &  
        \hspace{2cm} \times \frac{\left \| F_l(\ov{x}_{l-1};r_l,X_{r_l}^{r_{l-1},x_{l-1}}) -
                                F_l(\ov{x}_{l-1};r  ,X_r    ^{r_{l-1},x_{l-1}}) \right \|_p}{r_l-r} dr  \bigg ]\\
& =:& \varphi_l(\ov{x}_{l-1};s),
\tion
that means
\begin{equation}\label{eqn:upper_bound_U}
\| U_s^{r_{l-1},\ov{x}_{l-1}} \|_p
\le \varphi_l(\ov{x}_{l-1};s)
\end{equation}
with
\[     \noo \varphi_l(\ov{X}_{l-1};s) \rrm_p
   \le c_{(\ref{lemma:upper_bound_V})}  c_{(\ref{eqn:upper_bound_Al0})} 
       \bigg [  [r_l-s]+ \kappa_{p'} (1+\sqrt{T}) c_{B_{p,\infty}^\Theta} \int_s^{r_l}  
       (r_l-r)^{\frac{\theta_l}{2}-1} dr  \bigg ] \]
or
\begin{equation}\label{eqn:upper_bound_vph_l}
      \noo \varphi_l(\ov{X}_{l-1};s) \rrm_p
\le c_{(\ref{eqn:upper_bound_vph_l})} 
      \bigg [  [r_l-s]+ c_{B_{p,\infty}^\Theta} 
      \int_s^{r_l} (r_l-r)^{\frac{\theta_l}{2}-1} dr  \bigg ].
\end{equation}
Exploiting again Lemma \ref{lemma:upper_bound_V} also gives that
\equa
      \| V^{r_{l-1},\overline{x}_{r-1}}_s \|_p 
&\le& c_{(\ref{lemma:upper_bound_V})} \int_s^{r_l} \frac{\|A_l^0(\overline{x}_{l-1};r,X_r^{r_{l-1},x_{l-1}})\|_p}
                                              {\sqrt{r-s}} dr \\
&\le& c_{(\ref{lemma:upper_bound_V})} c_{(\ref{eqn:upper_bound_Al0})}
      \int_s^{r_l} \frac{ \psi_l(\ov{x}_{l-1};r)}
                                              {\sqrt{r-s}} dr
\tion
for $s\in [r_{l-1},r_l)\setminus \mathcal{N}_l(\overline{x}_{l-1})$, where $\mathcal{N}_l(\overline{x}_{l-1})$
has Lebesgue measure zero. Hence 
\equa
&   & \| U_s^{r_{l-1},\overline{x}_{l-1}} -U_t^{r_{l-1},\overline{x}_{l-1}} \|_p \\
& = & \bigg \| \int_s^t f^{\rm lin}(\overline{x}_{l-1};r,X_r^{r_{l-1},x_{l-1}},U_r^{r_{l-1},\overline{x}_{l-1}},V_r^{r_{l-1},\overline{x}_{l-1}})dr \\
&   &  \hspace*{15em}       - \int_s^t V_r^{r_{l-1},\overline{x}_{l-1}} dW_r^{r_{l-1}} \bigg \|_p \\
&\le& \left \| \int_s^t |A_l^0(\overline{x}_{l-1};r,X_r^{r_{l-1},x_{l-1}})| dr \right \|_p \\
&   &                + c_{(\ref{eqn:upper_bound_f_linear})}
      \left \| \int_s^t [|U_r^{r_{l-1},\overline{x}_{l-1}}| + |V_r^{r_{l-1},\overline{x}_{l-1}}|]  dr \right \|_p \\
&   & + a_p\noo \kla \int_s^t |V_r^{r_{l-1},\overline{x}_{l-1}}|^2 dr \mer^\frac{1}{2} \rrm_p \\
&\le& \left \| \int_s^t |A_l^0(\overline{x}_{l-1};r,X_r^{r_{l-1},x_{l-1}})| dr \right \|_p 
                + c_{(\ref{eqn:upper_bound_f_linear})}
      \int_s^t \|U_r^{r_{l-1},\overline{x}_{l-1}}\|_p dr \\
&   & + [c_{(\ref{eqn:upper_bound_f_linear})} \sqrt{t-s} + a_p]
      \kla \int_s^t \|V_r^{r_{l-1},\overline{x}_{l-1}}\|^2_p dr \mer^\frac{1}{2}\\
&\le& \left \| \int_s^t |A_l^0(\overline{x}_{l-1};r,X_r^{r_{l-1},x_{l-1}})| dr \right \|_p 
                + c_{(\ref{eqn:upper_bound_f_linear})}
      \int_s^t \varphi_l(\overline{x}_{l-1};r) dr \\
&   & + [c_{(\ref{eqn:upper_bound_f_linear})} \sqrt{t-s} + a_p]
         c_{(\ref{lemma:upper_bound_V})} c_{(\ref{eqn:upper_bound_Al0})}
         \kla \int_s^t 
         \bet \int_r^{r_l} \frac{ \psi_l(\ov{x}_{l-1};w)}{\sqrt{w-r}} dw \rag^2
         dr \mer^\frac{1}{2}.
\tion
Because $\P((X_{r_1},...,X_{r_{l-1}})\in E_l ) = 1$ we can use the stochastic flow property
and can bound $\|D_2(\overline{X}_{l-1})\|_p$ from above by the $L_p$-norms of the following
three expressions: Taking the $L_p$-norm of the last term gives
\equa
&   & \noo \kla \int_s^t 
         \bet \int_r^{r_l} \frac{ \psi_l(\ov{X}_{l-1};w)}{\sqrt{w-r}} dw \rag^2
         dr \mer^\frac{1}{2}\rrm_p \\
&\le& \kla \int_s^t 
                  \bet \int_r^{r_l}\!\! \frac{1 + \| \nabla_x F_l (\overline{X}_{l-1};w,X_w) \|_p
                 + \|      D^2 F_l (\overline{X}_{l-1};w,X_w) \|_p}{\sqrt{w-r}} dw\rag^2 \!\!
                         dr \mer^\frac{1}{2} \\
&\le&     \kla \int_s^t \bet \int_r^{r_l} \frac{dw}{\sqrt{w-r}} \rag^2
       dr \mer^\frac{1}{2} \\
&   & + \kappa_{p'} c_{B_{p,\infty}^\Theta} 
      \kla \int_s^t 
                  \bet \int_r^{r_l} \frac{ (r_l-w)^{\frac{\theta_l-1}{2}} + (r_l-w)^{\frac{\theta_l-2}{2}}
                 }{\sqrt{w-r}} dw\rag^2
                         dr \mer^\frac{1}{2} \\
&\le&     \kla \int_s^t \bet \int_r^{r_l} \frac{dw}{\sqrt{w-r}} \rag^2
       dr \mer^\frac{1}{2} \\
&   &  + \kappa_{p'} c_{B_{p,\infty}^\Theta} (1+\sqrt{T})
      \kla \int_s^t 
                  \bet \int_r^{r_l} \frac{  (r_l-w)^{\frac{\theta_l-2}{2}}
                 }{\sqrt{w-r}} dw\rag^2
                         dr \mer^\frac{1}{2} \\
&\le&   2 \sqrt{T} \sqrt{t-s} 
 + \kappa_{p'} c_{B_{p,\infty}^\Theta} (1+\sqrt{T})\gamma_l 
      \kla \int_s^t (r_l-r)^{\theta_l-1} dr \mer^\frac{1}{2} \\
\tion
with $\gamma_l:=\int_0^1 \frac{(1-t)^{\frac{\theta_l}{2}-1}}{\sqrt{t}} dt$.
For the next to the last term we obtain
\equa
      \noo \int_s^t \varphi_l(\overline{X}_{l-1};r) dr \rrm_p 
&\le& c_{(\ref{eqn:upper_bound_vph_l})} 
      \int_s^t \left [ (r_l-r) +  c_{B_{p,\infty}^\Theta} 
      \int_r^{r_l} 
         (r_l-w)^{\frac{\theta_l}{2}-1} dw \right ] dr  \\
&\le&  c_{(\ref{eqn:upper_bound_vph_l})}  \left [ T 
      + c_{B_{p,\infty}^\Theta} 
        \frac{2}{\theta_l} T^\frac{\theta_l}{2} \right ] (t-s).
\tion
Finally, we get by (\ref{eqn:upper_bound_psi_l}) and (\ref{eqn:smoothness_of_F_l})
that
\equa
&   &  \left \| \int_s^t  |A_l^0(\ov{X}_{l-1};r,X_r)| dr \right \|_p \\
&\le& c_{(\ref{eqn:upper_bound_Al0})} 
      \int_s^t  \|\psi_l(\ov{X}_{l-1};r)\|_p dr\\
&\le& c_{(\ref{eqn:upper_bound_Al0})} \left [ (t-s) + \sqrt{T} (1+\sqrt{T}) 
      \kappa_{p'} c_{B_{p,\infty}^\Theta}
      \kla \int_s^t (r_l-r)^{\theta_l-2} dr \mer^\frac{1}{2}
      \right ].
\tion

\underline{The term $D_3(\ov{X}_{l-1})$:} 
Let $r_{l-1} \le s < t < r_l$ and recall 
\[ Z_t^{r_{l-1},\overline{x}_{l-1}} = v_l(\overline{x}_{l-1};t,X_t^{r_{l-1},x_{l-1}})
                                       \sigma (t,X_t^{r_{l-1},x_{l-1}}).\]
From inequality (\ref{eqn:upper_bound_U}) we obtain 
\equa
&   & (t-s)^\frac{1}{2} \| Z_s^{r_{l-1},\overline{X}_{l-1}} \|_p \\
&\le& (t-s)^\frac{1}{2} \| \sigma \|_\infty \| v_l(\overline{X}_{l-1};s,X_s^{r_{l-1},X_{r_{l-1}}}) \|_p \\
&\le& (t-s)^\frac{1}{2} \| \sigma \|_\infty \left [ \| \nabla_x F_l(\overline{X}_{l-1};s,X_s^{r_{l-1},X_{r_{l-1}}}) \|_p 
                                            + \| U_s^{r_{l-1},\overline{X}_{l-1}} \|_p \right ] \\
&\le& (t-s)^\frac{ 1 }{2 }\|\sigma\|_{\infty} ( \kappa_{p'} c_{B_{p,\infty}^\Theta}
      (r_l-s)^\frac{\theta_l-1}{2} + \|\varphi_l(\ov{X}_{l-1},s)\|_p) \\
&\le& (t-s)^\frac{ 1 }{2 }\|\sigma\|_{\infty} \\
&   & \left ( \kappa_{p'} c_{B_{p,\infty}^\Theta}
        (r_l-s)^\frac{\theta_l-1}{2} + 
      c_{(\ref{eqn:upper_bound_vph_l})} 
      \bigg [  [r_l-s]+ c_{B_{p,\infty}^\Theta} 
      \int_s^{r_l} (r_l-r)^{\frac{\theta_l}{2}-1} dr  \bigg ] \right ) \\
&\le& c (t-s)^\frac{1}{2} [1+ (r_l-s)^\frac{\theta_l-1}{2}] \\
&\le& c \left [ (t-s)^\frac{1}{2} + \kla \int_s^t  (r_l-r)^{\theta_l-1} dr \mer^\frac{1}{2}
        \right ].
\tion

\underline{$(C3_l) \Longrightarrow (C6_l)$} 
Let
\[ t_k^{n,\theta_l} := r_{l-1} + (r_l-r_{l-1}) \kla 1 -  \kla 1 - \frac{k}{n} \mer^\frac{1}{\theta_l} \mer
   \sptext{1}{for}{1} 
  k=0,...,n\]
and
$S_{t_k^{n,\theta_l}}^n:= Y_{t_k^{n,\theta_l}}$. One obtains
for 
$t\in (t_{k-1}^{n,\theta_l},t_k^{n,\theta_l})
  \subseteq [r_{l-1},r_l]$
and an appropriate $\eta\in (0,1)$, that
\equa
&   & \| S_t^n - Y_t  \|_p \\
& = & \| (1-\eta)  Y_{t_{k-1}^{n,\theta_l}} + 
             \eta  Y_{t_k^{n,\theta_l}} - Y_t \|_p \\
&\le&   (1-\eta) \| Y_{t_{k-1}^{n,\theta_l}} - Y_t \|_p
      + \eta     \| Y_{t_k^{n,\theta_l}}     - Y_t \|_p \\
&\le& (1-\eta) c_3 \left ( \int_{t_{k-1}^{n,\theta_l}}^t  (r_l-r)^{\theta_l-1} dr 
           \right )^\frac{1}{2}
     + \eta c_3 \left ( \int_t^{t_k^{n,\theta_l}}  (r_l-r)^{\theta_l-1} dr 
         \right )^\frac{1}{2} \\
&\le&  c_3 \left ( \frac{1}{\theta_l}
      [ (r_l-t_{k-1}^{n,\theta_l})^{\theta_l}
       -(r_l-t_{k}^{n,\theta_l})^{\theta_l}]
          \right )^\frac{1}{2} \\
& = & c_3 \frac{(r_l-r_{l-1})^\frac{\theta_l}{2}}{\sqrt{ 
          \theta_l}} \frac{1}{\sqrt{n}}.
\tion

\underline{$(C6_l) \Longrightarrow (C4_l)$} We consider
\equa
&   & \left \|   Y_{\frac{r_l+t_{n-1}^{n,\theta_l}}{2}} 
               - S_{\frac{r_l+t_{n-1}^{n,\theta_l}}{2}}^n
      \right \|_p \\
& = & \left \|   Y_{\frac{r_l+t_{n-1}^{n,\theta_l}}{2}} - 
                 \frac{1}{2} \left [ 
                 S_{r_l}^n + S_{t_{n-1}^{n,\theta_l}}^n \right ]
      \right \|_p \\
&\ge& \left \|   Y_{\frac{r_l+t_{n-1}^{n,\theta_l}}{2}} - 
                 \frac{1}{2} \left [ 
                 Y_{r_l} + S_{t_{n-1}^{n,\theta_l}}^n \right ]
      \right \|_p -
      \frac{1}{2} \| Y_{r_l} - S_{r_l}^n \|_p 
\tion
so that
\[      \left \| Y_{r_l} - 2 Y_{\frac{r_l+t_{n-1}^{n,\theta_l}}{2}} 
                + S_{t_{n-1}^{n,\theta_l}}^n \right \|_p
   \le \frac{3c_6}{\sqrt{n}}. \]
But this means that
\[      \left \| Y_{r_l} - \E \kla Y_{r_l} | \ftn_{\frac{r_l+t_{n-1}^{n,\theta_l}}{2}} \mer
        \right \|_p 
   \le  \frac{6c_6}{\sqrt{n}}. \]
Because 
\[    r_l - \frac{r_l+t_{n-1}^{n,\theta_l}}{2} 
   = \frac{1}{2} (r_l-r_{l-1}) n^{-\frac{1}{\theta_l}} 
\]
we get that
\[   \left \| Y_{r_l} - \E \kla Y_{r_l} | \ftn_t \mer
        \right \|_p 
   \le  6c_6  \left (\frac{r_l-r_{l-1}}{2}\right )^{-\frac{\theta_l}{2}} (r_l-t)^{\frac{\theta_l}{2}}
   \sptext{1}{for}{1}
   t = \frac{r_l+t_{n-1}^{n,\theta_l}}{2}. \]
Using (\ref{eqn:lp:cond:expect}) proves  our assertion for 
$r_{l-1} + \frac{r_l-r_{l-1}}{2} \le t < r_l$. For the remaining 
$r_{l-1} \le t < r_{l-1} + \frac{r_l-r_{l-1}}{2}$ we can simply use
$\| Y_{r_l} - Y_{r_{l-1}}\|_p < \infty$.
\medskip

\underline{$(C7_l) \Longrightarrow (C4_l)$}
Let $t\in [r_{l-1},r_l)$. 
We use ($C7_l$) for $n=1$ so that $Y_t$ and $Y_{r_l}$ can be covered
by {\em one} ball with any radius bigger than $c_7(r_l-t)^\frac{\theta_l}{2}$. Taking the infimum
of these radii we get that 
$\| Y_{r_l}-Y_t\|_p \le 2 c_7(r_l-t)^\frac{\theta_l}{2}$ which implies that 
\[ \| Y_{r_l}- \E (Y_{r_l}|\ftn_t)\|_p \le 4 c_7(r_l-t)^\frac{\theta_l}{2}.\]
\underline{$(C3_l) \Longrightarrow (C7_l)$}
Fix $t\in [r_{l-1},r_l)$ and $n\ge 1$. Let $N\ge 1$ and choose $k\in \{ 1,...,N\}$ such that
\[ t \in       [t_{k-1}^{N,\theta_l},t_k^{N,\theta_l})
     \subseteq [r_{l-1},r_l). \]
For those time-nets we computed in $(C3) \Longrightarrow (C6)$ that
\[ \| Y_u-Y_v\|_p 
   \le c_3 \frac{(r_l-r_{l-1})^\frac{\theta_l}{2}}{\sqrt{\theta_l}} \frac{1}{\sqrt{N}} \] 
for $u,v\in [t_{k-1}^{N,\theta_l},t_k^{N,\theta_l}] \subseteq [r_{l-1},r_l]$. Now we choose 
$N\ge 1$ such that the cardinality of 
$\klae t_k^{N,\theta_l}:k=0,...,N\mere \cap [t_k^{N,\theta_l},r_l]$ is equal to $n$, i.e.
\[ n= 1 + N \left ( \frac{r_l -t_k^{N,\theta_l}}{r_l-r_{l-1}} \right )^{\theta_l}. \]
For $n\ge 2$ this implies that
\[ \frac{n}{2} \le n-1  
                =  \frac{N}{(r_l-r_{l-1})^{\theta_l}} (r_l-t_k^{N,\theta_l})^{\theta_l}
               \le \frac{N}{(r_l-r_{l-1})^{\theta_l}} (r_l-t)^{\theta_l}
\]
and
\[
      e_n ((Y_s)_{s\in [t,r_l]}|L_p)
 \le  c_3 \frac{(r_l-r_{l-1})^\frac{\theta_l}{2}}{\sqrt{ 
          \theta_l}} \frac{1}{\sqrt{N}}
 \le \frac{c_3}{\sqrt{ \theta_l}} 
      \frac{\sqrt{2 (r_l-t)^{\theta_l}}}{\sqrt{n}}. \]
The case $n=1$ implies that 
$t_{k-1}^{N,\theta_l}\le t < t_k^{N,\theta_l}=r_l$.
As in $(C3_l) \Longrightarrow (C4_l)$ we have 

\[     \| Y_{r_l} -     Y_s             \|_p 
  \le c_3 \sqrt{\frac{1}{\theta_l}} (r_l-s)^\frac{\theta_l}{2} 
  \le c_3 \sqrt{\frac{1}{\theta_l}} (r_l-t)^\frac{\theta_l}{2}
\]
for all $s\in [t,r_l]$ so that
\[     e_1 ((Y_s)_{s\in [t,r_l]}|L_p)
   \le c_3 \sqrt{\frac{1}{\theta_l}} (r_l-t)^\frac{\theta_l}{2}.
\]
\hfill$\Box$


\subsection{Proof of Theorem \ref{theorem:sufficient_I_new}}
\label{subsection:proof:sec:sufficient_conditions}
(a) We get, a.s., that
\begin{align*}
      X_s^\eta{}-X_s
& =    \int_0^s [b(r,X_r^\eta)-b(r,X_r)]dr \\
&    + \int_0^s [\sigma(r,X_r^\eta)-\sigma(r,X_r)]\sqrt{1-\eta(r)^2} dW_r \\
&    + \int_0^s \sigma(r,X_r^\eta) \eta(r) dB_r \\
&    - \int_0^s \sigma(r,X_r) (1-\sqrt{1-\eta(r)^2})         dW_r.
\end{align*}
Using the Burkholder-Davies-Gundy inequalities we estimate 
\[ e(s) := \E\sup_{0\le r\le s}|X^\eta_r-X_r|^p \]
by
\equa
      e(s)
&\le& 4^{p-1} \Big [ 
      T^{p-1} L_b^p\int_0^s e(r) dr+ a_p^p T^{p/2-1}L_\sigma^p\int_0^s e(r) dr \\
&   & + a_p^p \|\sigma\|_\infty^p \kla \int_0^s \eta(r)^2 dr \mer^\frac{p}{2}
      + a_p^p \|\sigma\|_\infty^p \kla \int_0^s (1-\sqrt{1-\eta(r)^2})^2 dr \mer^\frac{p}{2}
           \Big ],
\tion
where $L_b$ and $L_\sigma$ are the Lipschitz constants (with respect to $x$) of $b$ and $\sigma$, and
$a_p$ the constant from the Burkholder-Davis-Gundy inequality. 
Note that $1-\sqrt{1-\eta(r)^2}=\frac{\eta(r)^2}{1+\sqrt{1-\eta(r)^2}}\le |\eta(r)|$ 
using $|\eta(r)|\le 1$. Thus, applying Gronwall's lemma implies
\begin{equation}\label{eqn:estimate_difference_X} 
    \noo \sup_{0\le r\le s}|X^\eta_r-X_r|\rrm_p  
\le c_{(\ref{eqn:estimate_difference_X})} \kla \int_0^s \eta(r)^2 dr\mer^\frac{1}{2}
\end{equation}
where $c_{(\ref{eqn:estimate_difference_X})}>0$ depends at most on $(p,T,b,\sigma)$.
\bigskip

(b) We consider $Y^\eta-Y$ and $Z^\eta-Z$
and relate $(Y,Z)$ and $(Y^\eta,Z^\eta)$ to two BSDEs driven by the same 
Brownian motion $(W,B)$. This is the purpose of the construction below.
\smallskip

Let $\varphi=\chi_{[-1/2,1/2]}$ so that
\begin{equation}
  \label{eq:varphi}
     \sup_{\eta\in[-1,1]} \left (\frac{\varphi(\eta)}{\sqrt{1- \eta^2}}
                       +\frac{1-\varphi(\eta)}{|\eta|}\right ) = 2
\end{equation}
using the convention $\frac 00=0$. Thus, we can define the parameterized driver 
\[    f^\eta(t,\omega,y,z)
  := f\left (t,X^\eta_t(\omega),y,z^W\frac{\varphi(\eta(t))}
           {\sqrt{1- \eta(t)^2}}+z^B\frac{1-\varphi(\eta(t))}{\eta(t)}\right ) \]
where $z=(z^W,z^B)$ is $2d$-dimensional. In view of (\ref{eq:varphi}), the 
driver $f^\eta$ is Lipschitz with respect to $y$ and $z$. Thus, for 
any $\ftn^{W,B}_T$-measurable terminal condition $\widetilde{\xi} \in L_p$, there is an unique solution 
in $L_p$ in the filtration $\ftn^{W,B}$ to the BSDE
\[ \widetilde{Y}_t = \widetilde{\xi} + \int_t^T f^{\eta}(s,\widetilde{Y}_s,\widetilde{Z}_s) ds -
                                         \int_t^T \widetilde{Z}_s^W dW_s 
                                       - \int_t^T \widetilde{Z}_s^B dB_s  \]
because of \cite[Theorem 4.2]{Briand-et-al}.
\medskip

(c) For the driver $f^0$ (i.e. $\eta\equiv 0$) and terminal condition $\xi$
we have that
\[ (Y,[Z,0])\]
solves our BSDE.

(d) For the driver $f^\eta$ and the terminal condition $\xi^\eta$ we have that
\[  (Y^\eta,[Z^{\eta,W},Z^{\eta,B}]) \]
with
\[ Z^{\eta,W}_s=Z^\eta_s \sqrt{1-\eta(s)^2}
   \sptext{1}{and}{1}
   Z^{\eta,B}_s=Z^\eta_s \eta(s) \]
solves our BSDE because 
  \begin{align*}
    Y^\eta_t & =  \xi^\eta + \int_t^T f(s,X^\eta_s,Y^\eta_s, Z^\eta_s\varphi(\eta(s))+Z^\eta_s(1-\varphi(\eta(s)))) ds\\
& \qquad- \int_t^T Z^\eta_s \sqrt{1-\eta(s)^2}dW_s-\int_t^T Z^\eta_s \eta(s)dB_s\\
&=\xi^\eta + \int_t^T f^\eta(s,Y^\eta_s, [Z^{\eta,W}_s,Z^{\eta,B}_s])ds - \int_t^T Z^{\eta,W}_s dW_s-\int_t^T Z^{\eta,B}_s dB_s.
  \end{align*}

(e) To sum up, we have proved that $(Y,[Z,0])$ and $(Y^\eta,[Z^\eta_. \sqrt{1-\eta(.)^2},Z^\eta_. \eta(.)])$ solve
the BSDEs with data $(\xi,f^0)$ and $(\xi^\eta,f^\eta)$ in the filtration $(\ftn_t^{W,B})_{t\in [0,T]}$. 
Then, we are in a position 
to apply Lemma \ref{lemma:compare_briand_etal} (with $d$ replaced by $2d$) and get
\equa
&   &   \left \|\sup_{0\le t\le T}|Y^\eta_t-Y_t|\right \|_p
      + \left \|\kla \int_0^T \bet (Z^\eta_t\sqrt{1-\eta(t)^2}-Z_t,Z^\eta_t 
        \eta(t)) \rag^2 dt \mer^{1/2}
        \right \|_p\\
&\le& c_{(\ref{lemma:compare_briand_etal})} \left (\|\xi^\eta-\xi\|_p 
      + \noo \int_0^T |f^\eta(t,Y_t,[Z_t,0])-f^0(t,Y_t,[Z_t,0])|dt\rrm_p
      \right )\\
&\le& c_{(\ref{lemma:compare_briand_etal})} \bigg (
        \|\xi^\eta-\xi\|_p
      + L_f\noo\int_0^T |X^\eta_t-X_t|dt\rrm_p \\
&   & \hspace*{12em}
      + L_f\noo\int_0^T |Z_t|\left |\frac{\varphi(\eta(t))}{\sqrt{1- \eta(t)^2}}-1\right |dt\rrm_p
       \bigg )
\tion
where $c_{(\ref{lemma:compare_briand_etal})}$ (here and thereafter) is not identical with the constant 
$c$ in Lemma \ref{lemma:compare_briand_etal} but only refers to the fact that the inequality 
of Lemma \ref{lemma:compare_briand_etal} is used.
Now, since $\bet \frac{\varphi(\eta)}{\sqrt{1- \eta^2}}-1\rag \le c_{\varphi}|\eta|$ 
for some constant $c_{\varphi}>0$, we have 
\[     \int_0^T |Z_t|\bet \frac{\varphi(\eta(t))}{\sqrt{1- \eta(t)^2}}-1\rag dt
   \le c_\varphi \kla \int_0^T |Z_t|^2 dt \mer^{1/2} 
         \kla \int_0^T \eta(t)^2 dt \mer^{1/2}.\]
With the previous estimate on $X^\eta-X$ from (\ref{eqn:estimate_difference_X}) this leads to
\equa
&   & \noo \sup_{0\le t\le T}|Y^\eta_t-Y_t|\rrm_p
      +
          \noo \kla \int_0^T \bet (Z^\eta_t\sqrt{1-\eta(t)^2}-Z_t,Z^\eta_t \eta(t))\rag^2 dt
          \mer^\frac{1}{2}\rrm_p\\
&\le& c_{(\ref{lemma:compare_briand_etal})} \Bigg [ \|\xi^\eta-\xi\|_p \\
&   & +  L_f \left [ T c_{(\ref{eqn:estimate_difference_X})} 
      + c_\varphi \noo \kla \int_0^T |Z_t|^2 dt \mer^\frac{1}{2} \rrm_p \right ]
        \kla \int_0^T \eta(t)^2 dt \mer^\frac{1}{2} \Bigg ].
\tion
Applying Lemma \ref{lemma:compare_briand_etal} to $\xi^{(0)}=0$, $f_0=0$, $Y_s^{(0)}\equiv 0$, $Z_s^{(0)}\equiv 0$, 
$\xi^{(1)}=\xi$, $f_1(\omega;s,y,z):= f(s,X_s(\omega),y,z)$ and our solution $(Y,Z)$ we obtain
\[ \alpha_s(\omega)  =  |f(s,X_s(\omega),0,0)| 
                    \le K_f + L_f \sup_{0\le t \le T} |X_t(\omega)| \]
and
\[ \noo \kla \int_0^T |Z_t|^2 dt \mer^\frac{1}{2} \rrm_p
   \le c_{(\ref{lemma:compare_briand_etal})} [K_f+L_f + \|\xi\|_p]. \]
To complete the proof, it remains to use the inequality
\equa
      |(Z^\eta_t\sqrt{1-\eta(t)^2}-Z_t,Z^\eta_t \eta(t))|^2
& = & |Z^\eta_t|^2+|Z_t|^2-2 \sqrt{1-\eta(t)^2} \langle Z^\eta_t, Z_t \rangle \\
&\ge& \frac 12 |Z^\eta_t-Z_t|^2.
\tion
\hfill $\Box$


\subsection{Proof of Theorem \ref{thm:Lipschitz_functional}}
\label{sec:proof_of_Theorem_thm:Lipschitz_functional}

(a) In this step we assume that all $(x_1,...,x_L)$ (and similarly $(x'_1,...,x'_L)$) that appear have the property
that $x_1\in D_2$, $(x_1,x_2)\in D_3$, ..., $(x_1,...,x_{L-1})\in D_L$ where the
sets $D_2,...,D_L$ are taken from Proposition \ref{prop:Au}.
By backward induction we prove the following estimate regarding the terminal condition function $\Phi_{l}(x_1,...,x_{l})
   : = u_l\left (x_1,...,x_{l-1};r_{l},x_{l}\right )$ of the BSDE at time $r_{l}$:
   \begin{align}
&   \left | \Phi_{l}\left (\overline x _{l}\right )
            - \Phi_{l}\left (\overline x'_{l} \right ) \right | 
\le c_{l}         \sum_{i=1}^{l} \left [ \bet g_i(x_i) - g_i(x'_i) \rag 
       + \psi_i(\overline x_i; \overline x'_i) |x_i-x'_i| \right ].
\label{eq:lipschitz:terminal:condition}
\end{align}   
This is true for $l=L$ by our assumption. Assume now that \eqref{eq:lipschitz:terminal:condition} holds
for some $2\le l\leq L$ and let us prove the inequality for $l-1$. We have
\equa
&   & \left |   \Phi_{l-1} (x _1,...,x _{l-1})
              - \Phi_{l-1} (x'_1,...,x'_{l-1}) \right |\\
&\le& \big | u_l\left (x _1,...,x _{l-1};r_{l-1},x_{l-1}\right )            - u_l\left (x'_1,...,x'_{l-1};r_{l-1},x_{l-1}\right ) \big | \\
&   & + 
      \left | u_l\left (x'_1,...,x'_{l-1};r_{l-1},x _{l-1} \right ) 
            - u_l\left (x'_1,...         ,x'_{l-1};r_{l-1},x'_{l-1} \right ) \right  | \\
&\le& \big | u_l\left (x _1,...,x _{l-1};r_{l-1},x_{l-1}\right ) 
           - u_l\left (x'_1,...,x'_{l-1};r_{l-1},x_{l-1}\right ) \big | \\
&   & + \frac{\alpha_l}{\sqrt{r_l-r_{l-1}}} \\
&   &   (1+|x'_1|^{q_{l,1}}+\cdots+|x'_{l-1}|^{q_{l,l-1}}+|x_{l-1}|^{q_{l,l}} 
        + |x'_{l-1}|^{q_{l,l}}) 
              |x_{l-1}-x'_{l-1}|
\tion
where we used Proposition \ref{prop:Au}.
To estimate the remaining first term we use Lemma \ref{lemma:compare_briand_etal} and get that
\equa
&   & \left | u_l\left (x _1,...,x _{l-1};r_{l-1},x_{l-1} \right )
            - u_l\left (x'_1,...,x'_{l-1};r_{l-1},x_{l-1} \right ) \right | \\
&\le& c_{(\ref{lemma:compare_briand_etal})}
      \| u_l\left (x_1,..., x_{l-1};r_l,X_{r_l}^{r_{l-1},x_{l-1}}\right )      - u_l\left (x'_1,...,x'_{l-1};r_l,X_{r_l}^{r_{l-1},x_{l-1}}\right ) 
      \|_2 \\
&=& c_{(\ref{lemma:compare_briand_etal})}
      \| \Phi_l\left (x_1,..., x_{l-1},X_{r_l}^{r_{l-1},x_{l-1}}\right )      - \Phi_l\left (x'_1,...,x'_{l-1},X_{r_l}^{r_{l-1},x_{l-1}}\right )       \|_2 \\
&\le& c_{(\ref{lemma:compare_briand_etal})}c_{ \eqref{eq:lipschitz:terminal:condition} }
        \big( \sum_{i=1}^{l-1} \left [ \bet g_i(x_i) - g_i(x'_i) \rag 
       + \psi_i(x_1,...,x_i; x'_1,...,x'_i) |x_i-x'_i| \right ]\big).
\tion
(b) In the second step we verify the fractional smoothness, where we use (\ref{eqn:change_of_Y})
and therefore the  inequalities  from step (a). For $r_{l-1} \le s < r_l$, we have
\[  \| Y_{r_l} - \E (Y_{r_l}| \ftn_s ) \|_p 
  = \|      \Phi_l(X_{r_1},...,X_{r_l}) 
         - \E(\Phi_l(X_{r_1},...,X_{r_l})|\ftn_s)\|_p. \]
In particular, this expression 
depends on $x_0$, $b$, $\sigma$, $r_1,...,r_l$, $s$ and $\Phi_l$ but not on the
specific realization of the diffusion $X$. Hence we can assume the extended setting from Section \ref{sec:first_sufficient_condition}.
Using inequalities (\ref{eqn:conditional_expectation_vs_transition_density}) and 
the estimate \eqref{eq:lipschitz:terminal:condition} implies that
\equa
&   & \|      \Phi_l(X_{r_1},...,X_{r_l}) 
         - \E(\Phi_l(X_{r_1},...,X_{r_l})|\ftn_s)\|_p\\
&\le& \|      \Phi_l(X_{r_1},...,X_{r_l}) 
         -    \Phi_l(X_{r_1},...,X_{r_{l-1}},X_{r_l}^{\eta_{s,r_l}})\|_p\\
&\le& c_l \bigg \|   \bet g_l(X_{r_L}) - g_l(X_{r_l}^{\eta_{s,r_l}}) \rag \\
&   &    + \psi_l(X_{r_1},...,X_{r_l}; X_{r_1},...,X_{r_{l-1}},X_{r_l}^{\eta_{s,r_l}}) 
             |X_{r_l}-X_{r_l}^{\eta_{s,r_l}}|\bigg \|_p\\
&\le& c_l \| g_l(X_{r_l}) - g_l(X_{r_l}^{\eta_{s,r_l}}) \|_p \\
&   &  + c_l \| \psi_l(X_{r_1},...,X_{r_l}; X_{r_1},...,X_{r_{l-1}},X_{r_l}^{\eta_{s,r_l}})\|_{2p}
             \|X_{r_l}-X_{r_l}^{\eta_{s,r_l}}\|_{2p} \\
&\le& 2 c_l \| g_l(X_{r_l})- \E (g_l(X_{r_l})|\ftn_s) \|_p \\
&   &      + c_l \sup_{r_{l-1}\le u \le r_l} \| \psi_l(X_{r_1},...,X_{r_l};X_{r_1},...,X_{r_{l-1}},X_{r_l}^{\eta_{u,r_l}})\|_{2p}
       c_{(\ref{eqn:estimate_difference_X})} \sqrt{r_l-s}.
\tion
\hfill$\Box$


\section{Perspectives}

As natural steps, which could follow this paper, we see the investigation 
of more sufficient conditions for the fractional smoothness of a BSDE and 
the investigation of the limiting case as the number of points $r_1,...,r_L$ 
tends to infinity. In this connection the question, to what extend the generator 
might be path-dependent, is of interest as well. 
Moreover, the investigation of the above results in the context of other types of 
BSDEs (for example including reflection) and the development 
of numerical algorithms based on the discretizations proposed in this paper 
would be important.


\appendix
\section{Some lemmas about BSDEs}

We fix a complete probability space $(M,\Sigma,\Q)$,  
$0\le r < R \le T$ (the upper bound $T$ is used to bound some constants independently from $R$), 
$d\ge 1$ and a $d$-dimensional standard Brownian motion $B=(B_t)_{t\in [r,R]}$ 
with $B_r\equiv 0$. Furthermore, we assume that 
$(\gtn_t)_{t\in [r,R]}$ is the augmentation of the natural filtration of $B$.
The diffusion $(X_s)_{s\in [r,R]}$ is considered with respect to the same $\sigma$ and $b$ as used before,  restricted
to the corresponding time interval. Regarding the flow $(X_s^{t,x})_{s,t\in [r,R], x\in \R^d}$ and the
filtrations $(\gtn_s^t)_{s\in [t,R]}$ we use the same
convention as in Section \ref{sec:forward-backward}.

\begin{lemma}[$L_p$-stability of solutions of BSDEs]
\label{lemma:compare_briand_etal}
Let $2\le p < \infty$, 
$f_i:M\times [r,R]\times\R^k\times\R^{k\times d}\to \R^k$ be measurable with respect to
$Prog(M\times [r,R])\times \mathcal{B}(\R^k) \times \mathcal{B}(\R^{k\times d})$ with
$Prog(M\times [r,R])$ being the $\sigma$-algebra of progressively measurable 
subsets, and assume that, a.s.,
\[     Y_t^{(i)}
    =   \xi^{(i)} + \int_t^R f_i(s,Y_s^{(i)},Z_s^{(i)})ds 
      - \int_t^R Z_s^{(i)} d B_s
   \sptext{1}{for}{.3}
   i=0,1
   \sptext{.3}{and}{.3}
   r\le t \le R \]
with 
\[   \int_r^R |f_i(s,Y_s^{(i)},Z_s^{(i)})|ds
   + \sup_{r\le t \le R} |Y_t^{(i)}|
   + \kla \int_r^R|Z_s^{(i)}|^2 ds \mer^\frac{1}{2} \in L_p. \]
Let
\[ \alpha_s(\omega) := |   f_1(\omega;s,Y_s^{(0)}(\omega),Z_s^{(0)}(\omega))
                         - f_0(\omega;s,Y_s^{(0)}(\omega),Z_s^{(0)}(\omega))| \]
and suppose that there is a $L_{f_1}>0$ such that
\[ |f_1(\omega;s,u_1,v_1) - f_1(\omega;s,u_2,v_2)| \le L_{f_1} [ |u_1-u_2| + |v_1-v_2| ]. \]
Then there exists a $c_p>0$, depending on $p$ only, such that 
for $a\ge L_{f_1}  + L_{f_1}^2$ one has
\begin{multline*}
\E \left [ \sup_{t\in [r,R]} e^{ap(t-r)}|\Delta Y_t|^p + \kla \int_r^R e^{2a(s-r)}|\Delta Z_s|^2 ds \mer^\frac{p}{2}
   \right ] \\
\le c_p^p \E \left [ e^{ap(R-r)}|\Delta \xi|^p + \kla \int_r^R e^{a(s-r)} \alpha_s ds \mer^p
\right ].
\end{multline*}
\end{lemma}

\begin{proof}
The result is a direct consequence of \cite[Proposition 3.2]{Briand-et-al}.
For $\Delta Y_t := Y_t^1 - Y_t^0$, $\Delta Z_t := Z_t^1 - Z_t^0$ 
and $\Delta \xi :=\xi^{(1)}-\xi^{(0)}$ we get that
\equa
      \Delta Y_t 
& = & \Delta \xi + \int_t^R \widehat{f}(s,\Delta Y_s,\Delta Z_s) ds 
                             - \int_t^R \Delta Z_s d B_s
\tion
with
$    \widehat{f}(s,\Delta y,\Delta z) 
 :=  f_1(s,\Delta y + Y_s^{(0)}, \Delta z + Z_s^{(0)})- f_0(s,Y_s^{(0)},Z_s^{(0)})$ and
 \equa
&   &
      |\widehat{f}(\omega;s,\Delta y,\Delta z)| \\
& = & |   f_1(\omega;s,\Delta y + Y_s^{(0)}(\omega), \Delta z + Z_s^{(0)}(\omega))
        - f_0(\omega;s,Y_s^{(0)}(\omega),Z_s^{(0)}(\omega))| \\
&\le& |   f_1(\omega;s,Y_s^{(0)}(\omega),Z_s^{(0)}(\omega))
        - f_0(\omega;s,Y_s^{(0)}(\omega),Z_s^{(0)}(\omega))| \\
&   & + | f_1(\omega;s,\Delta y + Y_s^{(0)}(\omega), \Delta z + Z_s^{(0)}(\omega))
        - f_1(\omega;s,Y_s^{(0)}(\omega),Z_s^{(0)}(\omega))| \\
&\le& \alpha_s(\omega) + L_{f_1} [|\Delta y| + |\Delta z|].
\tion
Applying \cite[Proposition 3.2]{Briand-et-al} implies the assertion.
\end{proof}

The following lemma shows that \cite[Theorem 3.2]{ZhangJ2} transfers to our path dependent setting as expected.
The proof is presumably only included in this preprint version
for the convenience of the reader as it is a copy of that one in
\cite{ZhangJ2} (see also \cite[Section 5]{Hu-Ma:2004}).

\begin{lemma}[Representation of a BSDE parameterized by a parameter $y\!\in\!\! \R^K$]
\label{lemma:validity_of_assumptions}
Assume that $(A_{b,\sigma})$ and $(A_f)$ are satisfied, that $K,d\ge 1$ and that 
$H: \R^K \times \R^d \to \R$ is Borel-measurable with
\[ |H(y;x)| \le \alpha (1 + |y|^\gamma + |x|^\beta)=: \psi(y,x) \]
for some $\alpha,\beta,\gamma \in [1,\infty)$. Then there exists a Borel set $F\subseteq \R^K$
such that $F^c$ is of Lebesgue measure zero and such that for
\[ G(y;x) := \chi_F(y) H(y;x) \]
and
\[   U(y;t,x):= \left \{ \begin{array}{lcl}
                Y_t^{y;t,x} \mbox{ a.s.} &:& r \le t < R \\
                    G(y;x)  &:& t = R 
                          \end{array} \right ., \]   
where $(Y_s^{y;t,x})_{s\in [t,R]}$ is the $Y$-component 
of the BSDE with respect to the forward diffusion $(X_s^{t,x})_{s\in [t,R]}$, 
the terminal condition $G(y;X_R^{t,x})$ with terminal time $R\in(0,T]$, and the generator $f$, the following 
assertions are satisfied:
\begin{enumerate}[{\rm (i)}]
\item For fixed $y\in \R^K$ we have that $U(y;\cdot,\cdot)\in C^{0,1}([r,R)\times \R^d)$.
\item The functions $         U:\R^K\times [r,R]\times \R^d \to \R  $ and
                    $\nabla_x U:\R^K\times [r,R)\times \R^d \to \R^{1\times d}$  
      are measurable.
\item There exists a constant $c>0$ depending at most on $(b,\sigma,T,\alpha,\gamma,\beta,K_f,L_f)$
      such that 
      \begin{enumerate}[{\rm (a)}]
      \item $|U(y;t,x)| \le c \psi(y;x)$ for $(y,t,x)\in \R^K\times [r,R] \times \R^d$,
      \item $|\nabla_x U(y;t,x)| \le c \frac{\psi(y;x)}{\sqrt{R-t}}$ for 
            $(y,t,x)\in \R^K\times [r,R) \times \R^d$.
      \end{enumerate}
\item For any $y\in \R^K$, the solution of the BSDE with the terminal condition $G(y;X_R^{r,x})$, generator $f$, 
      and forward diffusion $(X_s^{r,x})_{s\in [r,R]}$
      can be represented as 
      \begin{enumerate}[{\rm (a)}]
      \item $Y_t^{y;r,x}= U(y;t,X_t^{r,x})$ on $[r,R]$,
      \item $Z_t^{y;r,x}=\nabla_x U(y;t,X_t^{r,x})\sigma(t,X_t^{r,x})$ on $[r,R)$.
      \end{enumerate}
\end{enumerate}
\end{lemma}
\bigskip

\begin{proof}
We find 
$H_n\in C^\infty_0(\R^K\times \R^d)$, $n\ge 1$, such that 
\[ \lim_n H_n = H \quad \lambda_{K+d} \mbox{-a.e.} 
   \sptext{1}{and}{1}
  |H_n(y;x)| \le 2 \psi (y,x). \]
Hence there is a  Borel set $F\subseteq \R^K$
such that $F^c$ is of Lebesgue measure zero and such that
\[ \lim_n G_n(y;\cdot) = G(y;\cdot) \quad \lambda_{d} \mbox{-a.e.} \]
for all $y\in \R^K$ with
\[  G_n(y;x) := \chi_F(y) H_n (y;x) 
    \sptext{1}{and}{1}
    G(y;x) := \chi_F(y) H(y;x).\]
Let $U^n$ be defined as $U$ with $G$ replaced by $G_n$.
Applying \cite[Theorems 3.1 and 4.2]{Ma-Zhang:2002} gives that
\[ (U^n(y;s,X_s^{t,x}), \nabla_x U^n(y;s,X_s^{t,x}) \sigma(s,X_s^{t,x}))_{s\in [t,R]} \]
solves our BSDE on the interval $[t,R]$, that $U^n(y;\cdot,\cdot)\in C^{0,1}([r,R]\times \R^d)$ 
and that
\begin{multline*}
     \nabla_x U^n(y;t,x) 
   = \\ \E \left [   G_n(y; X_R^{t,x}) N_R^{t,1,(t,x)} \!
                + \! \int_t^R \!\! f(s,X_s^{t,x},Y_s^{y,n;t,x},Z_s^{y,n;t,x}) N^{t,1,(t,x)}_s ds
               \right ].
\end{multline*}

\underline{Properties of the function $U^n$}
\smallskip

{\bf (a)} 
To estimate $U^n(y;t,x)$ we let $U^n_0(y;t,x)$ be the corresponding solution
with the zero generator and denote by $(Y_{s,0}^{y,n;t,x},Z_{s,0}^{y,n;t,x})_{s\in [t,R]}$
the corresponding solution to our BSDE. By Lemma \ref{lemma:compare_briand_etal} we get that
\equa
&   & |U^n(y;t,x) - U^n_0(y;t,x)|\\
&\le& c_{(\ref{lemma:compare_briand_etal})}
      \noo  \int_t^R |f(s,X_s^{t,x},Y_{s,0}^{y,n;t,x},Z_{s,0}^{y,n;t,x})|ds \rrm_2 \\
&\le& c_{(\ref{lemma:compare_briand_etal})} K_f R + c_{(\ref{lemma:compare_briand_etal})} L_f 
      \noo\int_t^R [|X_s^{t,x}|+|Y_{s,0}^{y,n;t,x}|+|Z_{s,0}^{y,n;t,x}|]ds \rrm_2 \\
&\le& c_{(\ref{lemma:compare_briand_etal})} K_f R + c_{(\ref{lemma:compare_briand_etal})} L_f \times \\
&   & \hspace*{-1.1em}
      \noo \int_t^R [  |X_s^{t,x}|+|\E(G_n(y;X_R^{t,x})|\gtn_s^t)|
                        + \| \sigma \|_\infty |\E(G_n(y;X_R^{t,x})N_R^{s,1,(t,x)}|\gtn_s^t)]ds \rrm_2 \\
&\le& c_{(\ref{lemma:compare_briand_etal})} K_f R + c_{(\ref{lemma:compare_briand_etal})} L_f \times \\
&   & \hspace*{-1em}
      \noo \int_t^R [  |X_s^{t,x}| + 2 |\E(\psi(y;X_R^{t,x})^2|\gtn_s^t)|^\frac{1}{2}
                           (1+ \| \sigma \|_\infty|\E((N_R^{s,1,(t,x)})^2|\gtn_s^t)|^\frac{1}{2})
                        ]ds \rrm_2.
\tion
Now we use that
\[     \E((N_R^{s,1,(t,x)})^2|\gtn_s^t)
   \le \frac{\kappa_2^2}{R-s}
   \sptext{1}{and}{1}
   |U^n_0(y;t,x)| \le  2 \E \psi(y;X_R^{t,x}) \]
and
\[ \| \psi(y;X_s^{t,x}) \|_2 \le \beta_2 \psi(y;x)
   \sptext{1}{and}{1}
   \|X_s^{t,x}\|_2 \le \alpha_2 [1+|x|] \]
to get 
\equa
&   & |U^n(y;t,x)| \\
&\le& 2 \E \psi(y;X_R^{t,x}) + c_{(\ref{lemma:compare_briand_etal})} K_f R
       + c_{(\ref{lemma:compare_briand_etal})} L_f \times \\
&   & \noo \int_t^R [  |X_s^{t,x}|+ 2 |\E(\psi(y;X_R^{t,x})^2|\gtn_s^t)|^\frac{1}{2}
                           (1+ \| \sigma \|_\infty\kappa_2 (R-s)^{-1/2})
                        ]ds \rrm_2 \\                      
&\le& 2 \E \psi(y;X_R^{t,x}) + c_{(\ref{lemma:compare_briand_etal})} K_f R
       + c_{(\ref{lemma:compare_briand_etal})} L_f \times \\
&   & \left [   \int_t^R  \|X_s^{t,x}\|_2 ds  
              + 2 \int_t^R [\| \psi(y;X_R^{t,x})\|_2
                           (1+ \| \sigma \|_\infty \kappa_2 (R-s)^{-1/2})
                       ds \right ] \\                      
&\le& \| \psi(y;X_R^{t,x})\|_2 
      \left [ 2 + 2 c_{(\ref{lemma:compare_briand_etal})} L_f [R + 2 \| \sigma \|_\infty \kappa_2 R^{1/2}] \right ] \\
&   &  + c_{(\ref{lemma:compare_briand_etal})} L_f \int_t^R  \|X_s^{t,x}\|_2 ds 
       + c_{(\ref{lemma:compare_briand_etal})} K_f R \\
&\le& \beta_2 \psi(y;x)
      \left [ 2 + 2 c_{(\ref{lemma:compare_briand_etal})} L_f [R + 2 \| \sigma \|_\infty \kappa_2 R^{1/2}] \right ] \\
&   &  + c_{(\ref{lemma:compare_briand_etal})} L_f \int_t^R  \alpha_2 [1+|x|]  ds 
       + c_{(\ref{lemma:compare_briand_etal})} K_f R
\tion
so that, for some $c_{(\ref{eqn:upper_bound_Un})} \ge 1$,
\begin{equation}\label{eqn:upper_bound_Un}
|U^n(y;t,x)| \le  c_{(\ref{eqn:upper_bound_Un})} \psi(y;x). 
\end{equation}

{\bf (b)} According to \cite[Theorem 3.1, Corollary 3.2]{Ma-Zhang:2002} the gradient $\nabla_x U^n(y;t,x)$ exists 
and is bounded by a constant that might depend on $n$ and $y$. For $r\le t \le \rho \le R$ and $y \in \R^K$ we define 
\equa
A_t^\rho &:= & \sqrt{\rho-t} \sup_{x\in  \R^d} \frac{|\nabla_x U^n(y;t,x)|}{\psi(y;x)}, \\
B_t^\rho &:= & \sup_{s\in [t,\rho]} A_s^\rho.
\tion
Although $A^\rho$ and $B^\rho$ might depend on $n$ and $y$, we do not indicate this for the purpose of notational simplicity.
Using \cite[Theorem 4.2]{Ma-Zhang:2002} yields 
\equa
&   & |\nabla_x U^n(y;t,x)| \\
& = & \left | \E \left [ U^n(y;\rho, X_\rho^{t,x}) N_\rho^{t,1,(t,x)} + \! \int_t^\rho f(s,X_s^{t,x},Y_s^{y,n;t,x},Z_s^{y,n;t,x}) 
      N^{t,1,(t,x)}_s ds \right ] \right | \\
&\le&  c_{(\ref{eqn:upper_bound_Un})} \kappa_2 \frac{\|\psi(y;X_\rho^{t,x})\|_2}{\sqrt{\rho-t}}  + 
      \kappa_2 \int_t^\rho \frac{\| f(s,X_s^{t,x},Y_s^{y,n;t,x},Z_s^{y,n;t,x}) \|_2}{\sqrt{s-t}} ds \\
&\le& c_{(\ref{eqn:upper_bound_Un})} \kappa_2 \frac{\|\psi(y;X_\rho^{t,x})\|_2}{\sqrt{\rho-t}} 
      + \kappa_2 \times \\
&   & \hspace*{-3em} \int_t^\rho \frac{K_f + L_f  [ \|X_s^{t,x}\|_2 + c_{(\ref{eqn:upper_bound_Un})} \| \psi(y;X_s^{t,x}) \|_2
      + \| \nabla_x U^n(y;s,X_s^{t,x}) \sigma(s,X_s^{t,x}) \|_2]}{\sqrt{s-t}} ds \\
&\le& c_{(\ref{eqn:upper_bound_Un})}  \kappa_2 \frac{\|\psi(y;X_\rho^{t,x})\|_2}{\sqrt{\rho-t}} 
      + \kappa_2 \times \\
&   &  \hspace*{-2em}\int_t^\rho \frac{K_f + L_f  [ \|X_s^{t,x}\|_2 + c_{(\ref{eqn:upper_bound_Un})}  \| \psi(y;X_s^{t,x}) \|_2
      + \|\sigma\|_\infty \| \nabla_x U^n(y;s,X_s^{t,x}) \|_2]}{\sqrt{s-t}} ds \\  
&\le& c_{(\ref{eqn:upper_bound_Un})} \kappa_2 \frac{\|\psi(y;X_\rho^{t,x})\|_2}{\sqrt{\rho-t}} 
      + \\
&   & \hspace*{-2em} \kappa_2 \int_t^\rho\frac{K_f + L_f  \left [ \|X_s^{t,x}\|_2 + c_{(\ref{eqn:upper_bound_Un})}  \| \psi(y;X_s^{t,x}) \|_2
      + \|\sigma\|_\infty \left \| A_s^\rho \frac{\psi(y;X_s^{t,x})}{\sqrt{\rho-s}} \right \|_2\right ]}
        {\sqrt{s-t}} ds \\  
&\le&  c_{(\ref{eqn:upper_bound_Un})} \kappa_2 \frac{\|\psi(y;X_\rho^{t,x})\|_2}{\sqrt{\rho-t}} \\
&   &  + \kappa_2 2 \sqrt{\rho-r} [K_f+L_f+ L_f c_{(\ref{eqn:upper_bound_Un})} ] \sup_{s\in [t,\rho]} [ 1 + \|X_s^{t,x}\|_2 + \| \psi(y;X_s^{t,x}) \|_2 ] \\
&   & + \kappa_2 \|\sigma\|_\infty B_t^\rho  
      \sup_{s\in [t,\rho]} \| \psi(y;X_s^{t,x}) \|_2  B\kla\frac{1}{2},\frac{1}{2}\mer\\
&\le&  c_{(\ref{eqn:upper_bound_Un})}  \kappa_2 \frac{\beta_2 \psi(y;x)}{\sqrt{\rho-t}} + 
      \kappa_2 2 \sqrt{\rho-r} [K_f+L_f+ L_f c_{(\ref{eqn:upper_bound_Un})} ] [ 1 + \alpha_2 [1+|x|] \\
&   & + \beta_2 \psi(y;x) ] + \kappa_2 \|\sigma\|_\infty B_t^\rho  \beta_2 \psi(y;x)  B\kla \frac{1}{2},\frac{1}{2}\mer\\
&\le&  A_\psi \psi(y;x) \kla \frac{1}{\sqrt{\rho-t}} + B_t^\rho \mer
\tion
where $A_\psi>0$ depends at most on $(b,\sigma,T,\alpha,\beta,\gamma,K_f,L_f)$.
Consequently,
\[     A_t^\rho  
   \le A_\psi \kla 1 + \sqrt{\rho-t} B_t^\rho \mer \]
and
\[     B_t^\rho  
   \le A_\psi \kla 1 + \sqrt{\rho-t} B_t^\rho \mer. \]
In case of $|\rho-t|\le (2  A_\psi)^{-2}$ this gives
$B_t^\rho \le 2 A_\psi$ and
\[ |\nabla_x U^n(y;t,x)| \le 2 A_\psi \frac{\psi(y;x)}{\sqrt{\rho-t}}. \]
Moreover, in case of $\frac{1}{4} (2  A_\psi)^{-2} \le |\rho-t|\le (2  A_\psi)^{-2}$ we also get that
\[ |\nabla_x U^n(y;t,x)| \le 2 A_\psi  \frac{\psi(y;t,x)}{\sqrt{\rho-t}}
                         \le 8 A_\psi^2 \psi(y;x).\]
The latter inequality means that 
\[ |\nabla_x U^n(y;t,x)| \le 8 A_\psi^2 \psi(y;x)\]
whenever $r\le t \le R$ and $|R-t|\ge \frac{1}{4} (2  A_\psi)^{-2}$. On the other side,
\[ |\nabla_x U^n(y;t,x)| \le 2 A_\psi \frac{\psi(y;x)}{\sqrt{R-t}} \]
for $r\le t \le R$ and $|R-t|\le (2  A_\psi)^{-2}$. Combining both estimates yields to
\begin{equation}\label{eqn:upper_bound_nabla_Un}
|\nabla_x U^n(y;t,x)| \le c_{(\ref{eqn:upper_bound_nabla_Un})} \frac{\psi(y;x)}{\sqrt{R-t}}
\end{equation}
for all $t\in [r,R]$.
\smallskip

{\bf (c)} We show that $U^n(y;t,x)$ is measurable as function on $\R^K\times [r,R] \times \R^d$.
Let $y,y'\in F$. Then it follows by Lemma \ref{lemma:compare_briand_etal} that
\equa
      |U^n(y;t,x) - U^n(y';t,x)|
&\le& c_{(\ref{lemma:compare_briand_etal})}
      \noo G_n(y;X_R^{(t,x)}) -  G_n(y';X_R^{(t,x)}) \rrm_2 \\
&\le& c_{(\ref{lemma:compare_briand_etal})} Lip(H^n) |y-y'|
\tion
and
\equa
&   & |U^n(y;t,x) - U^n(y';t',x')|\\
&\le& |U^n(y;t,x) - U^n(y';t,x)| + |U^n(y';t,x) - U^n(y';t',x')| \\
&\le& c_{(\ref{lemma:compare_briand_etal})} Lip(H_n) |y-y'|
      + |U^n(y';t,x) - U^n(y';t',x')|.
\tion
Hence 
$(U^n)^{-1}(B) \cap (F\times [r,R]\times \R^d) \in {\mathcal{B}}(\R^K\times [r,R]\times \R^d)$
for all open sets $B\subseteq \R$. On $F^c$ we have 
\[   (U^n)^{-1}(B) \cap (F^c\times [r,R]\times \R^d)
   = F^c \times \overline{U}^{-1}(B) \in {\mathcal{B}}(\R^K\times [r,R]\times \R^d) \]
where $\overline{U}(t,x)$ is the functional for the $Y$ process with zero 
terminal condition. Consequently, $U^n$ is measurable.
\medskip

\underline{Properties of the function $U$} 
\smallskip

{\bf (d)} Let $D$ be the product of $[r,b]\subseteq [r,R)$ where $b\in (r,R)$
and a compact subset of $\R^d$. For $(t,x)\in D$ Lemma \ref{lemma:compare_briand_etal}
and Proposition \ref{proposition:transition_density} yield
\equa
&   & |U(y;t,x)-U^n(y;t,x)|^2\\
&\le& c_{(\ref{lemma:compare_briand_etal})}^2 \|G(y;X_R^{t,x})-G_n(y;X_R^{t,x})\|_2^2 \\
& = & c_{(\ref{lemma:compare_briand_etal})}^2 
      \int_{\R^d} \Gamma(t,x;R,\xi) |G(y;\xi)-G_n(y;\xi)|^2 d\xi \\
&\le& c_{(\ref{lemma:compare_briand_etal})}^2 
      \int_{\R^d} c_{(\ref{proposition:transition_density})}
      \gamma_{R-t}^d\kla \frac{x-\xi}{c_{(\ref{proposition:transition_density})}} \mer
      |G(y;\xi)-G_n(y;\xi)|^2 d\xi \\
&\le& c_{(\ref{lemma:compare_briand_etal})}^2 
      \int_{\R^d} \frac{c_{(\ref{proposition:transition_density})}}
                       {(2\pi(R-b))^{\frac{d}{2}}}
      e^{-\frac{1}{c_{(\ref{proposition:transition_density})}^2}\frac{|x-\xi|^2}{(R-r)}} 
      |G(y;\xi)-G_n(y;\xi)|^2 d\xi \\
&\le& c_{(\ref{lemma:compare_briand_etal})}^2 
      \int_{\R^d} \frac{c_{(\ref{proposition:transition_density})}}
                       {(2\pi(R-b))^{\frac{d}{2}}}
      e^{\frac{1}{c_{(\ref{proposition:transition_density})}^2}
          \frac{- (|\xi|^2/2) + |x|^2}{(R-r)}} 
      |G(y;\xi)-G_n(y;\xi)|^2 d\xi.
\tion
This implies that for all fixed parameters $y\in \R^K$ there is a uniform
convergence on $D$ of $U^n$ towards $U$, so that 
$U(y;\cdot,\cdot)$ is continuous on $[r,R)\times \R^d$.
Moreover, as limit of measurable functions $U^n$, the
function $U:\R^K\times [r,R)\times \R^d\to \R$ is measurable as well. Because 
$U(y;R,x) = G(y;x)$ the function $U:\R^K\times [r,R]\times \R^d\to \R$ is
measurable.
\smallskip

{\bf (e)} Now we get
\equa
&   & \| Y_t^{y;r,x} - U(y;t,X_t^{r,x})\|_2 \\
&\le& \| Y_t^{y;r,x} - U^n(y;t,X_t^{r,x})\|_2 + \| U^n(y;t,X_t^{r,x}) - U(y;t,X_t^{r,x})\|_2\\
&\le& c_{(\ref{lemma:compare_briand_etal})} \|G(y;X_R^{r,x})-G_n(y;X_R^{r,x})\|_2 +
      \| U^n(y;t,X_t^{r,x}) - U(y;t,X_t^{r,x})\|_2
\tion
where we have used Lemma \ref{lemma:compare_briand_etal}.
Using dominated convergence both terms converge to zero as $n\to \infty$, because
\[ |U(y;t,x)| \le c_{(\ref{eqn:upper_bound_Un})}\psi(y;x)\]
as a consequence of (\ref{eqn:upper_bound_Un}) and step (d). Consequently,
\[ Y_t^{y;r,x} = U(y;t,X_t^{r,x}) \quad a.s. 
   \sptext{1}{for all}{1}
   t\in [r,R]. \]
\underline{The function $\nabla_x U(y;t,x)$} 
\smallskip

{\bf (f)} By Lemma \ref{lemma:compare_briand_etal} we know that
\equa
&   & \lim_n \int_t^R \noo \nabla_x U^n(y;s,X^{t,x}_s) \sigma(s,X_s^{t,x}) - 
      Z_s^{y;t,x} \rrm_2^2 ds \\
& = & \lim_n \int_t^R \noo Z_s^{y,n;t,x} - Z_s^{y;t,x} \rrm_2^2 ds = 0.
\tion
Let 
\[    V(y;t,x) 
   := \E \left [ G(y; X_R^{t,x}) N_R^{t,1,(t,x)} + \int_t^R f(s,X_s^{t,x},Y_s^{y;t,x},Z_s^{y;t,x}) 
                 N^{t,1,(t,x)}_s ds
               \right ]. \]
By dominated convergence we have that
\equa
&   & \lim_n \nabla_x U^n(y;t,x) \\
& = & \lim_n  \E \bigg  [   G_n(y; X_R^{t,x}) N_R^{t,1,(t,x)} \\
&   & + \! \int_t^R \!\! 
         f\big (s,X_s^{t,x},U^n(y;s,X_s^{t,x}),\nabla_x U^n (y;s,X_s^{t,x})\sigma(s,X_s^{t,x})
          \big )
            N^{t,1,(t,x)}_s ds
               \bigg ]\\
& = & \E \bigg  [   G(y; X_R^{t,x}) N_R^{t,1,(t,x)}  + \! \int_t^R \!\! 
         f\big (s,X_s^{t,x},Y_s^{y;t,x},Z_s^{y;t,x}
          \big )
            N^{t,1,(t,x)}_s ds
               \bigg ]\\
& = &  V(y;t,x) 
\tion
for all $(t,x)\in [r,R)\times \R^d$, which also implies
\[ |V(y;t,x)| \le c_{(\ref{eqn:upper_bound_nabla_Un})}  \frac{\psi(y;x)}{\sqrt{R-t}} \] 
by (\ref{eqn:upper_bound_nabla_Un}). Consequently,
\[    \lim_n \int_r^{R-\delta} \noo \nabla_x U^n(y;t,X_t^{r,x}) \sigma(t,X_t^{r,x}) - 
                                    V  (y;t,X_t^{r,x}) \sigma(t,X_t^{r,x})  \rrm_2^2 dt = 0 \]
for all $\delta \in (0,R-r)$ and 
\[ Z_t^{y;r,x} =   V  (y;t,X^{r,x}_t) \sigma(t,X_t^{r,x})  \quad a.s. \]
for almost every $t\in [r,R)$. So we can re-define
\[  Z_t^{y;r,x} := V  (y;t,X^{r,x}_t) \sigma(t,X_t^{r,x}). \]
{\bf (g)} Next we show that 
\[    V(y;t,x) 
   = \E \left [ G(y; X_R^{t,x}) N_R^{t,1,(t,x)} + \int_t^R f(s,X_s^{t,x},Y_s^{y;t,x},Z_s^{y;t,x}) N^{t,1,(t,x)}_s ds
               \right ] \]
is continuous in $(t,x)$ on $[r,R)\times \R^d$. 
For the first term this follows
from the classical theory from the properties of the transition density of 
$X$ because  
\[ \E \left [ G(y; X_R^{t,x}) N_R^{t,1,(t,x)}\right]=\int_{\R^d}G(y;w) \nabla_x \Gamma(t,x;R,w)dw
   \] 
and the continuity in $(t,x)$ follows from Proposition \ref{proposition:transition_density}.
So it remains to show that
\[ (t,x) \to  \E \int_t^R f(s,X_s^{t,x},Y_s^{y;t,x},Z_s^{y;t,x}) N^{t,1,(t,x)}_s ds \]
is continuous in $(t,x)$ on each  $D$ which is the product of $[r,b]\subseteq [r,R)$ 
and a compact subset of $\R^d$. Take a sequence $(t_n,x_n)\to (t,x)$ from $D$. We write
\equa
&   & \E \int_t^R f(s,X_s^{t,x},Y_s^{y;t,x},Z_s^{y;t,x}) N^{t,1,(t,x)}_s ds \\
& = &    \int_{(r,R)} \frac{\chi_{(t,R)}(s)}{\sqrt{R-s}\sqrt{s-t}} \\
&   &                  \E \big [ [\sqrt{R-s} f(s,X_s^{t,x},Y_s^{y;t,x},Z_s^{y;t,x})] 
                        [\sqrt{s-t} N^{t,1,(t,x)}_s] \big ] ds \\
& = & \int_{(r,R)} \varphi_t(s) \psi_{t,x}(s) ds  
\tion
with
\equa
\varphi_t(s)  & := & \frac{\chi_{(t,R)}(s)}{\sqrt{R-s}\sqrt{s-t}},\\
\psi_{t,x}(s) & := & \chi_{(t,R)}(s) \E \big [ [\sqrt{R-s} f(s,X_s^{t,x},Y_s^{y;t,x},Z_s^{y;t,x})] 
                        [\sqrt{s-t} N^{t,1,(t,x)}_s] \big ].
\tion
The family $(\varphi_t)_{t\in [r,b]}$ is uniformly integrable for $b\in [r,R)$. 
The boundedness of $(\psi_{t,x})_{(t,x)\in D}$ follows from 
\equa
      |\psi_{t,x}(s)|
&\le& \sqrt{R-s}  \| f(s,X_s^{t,x},Y_s^{y;t,x},Z_s^{y;t,x}) \|_2
      \sqrt{s-t}  \| N^{t,1,(t,x)}_s \|_2 \\
&\le& \sqrt{R-s} \bigg [ K_f + L_f (\| X_s^{t,x}\|_2 + \| U(y;s,X_s^{t,x})\|_2 \\
&   &  + \|\sigma\|_\infty \| V(y;s,X_s^{t,x})\|_2)\bigg ] \kappa_2
\tion
and the previous estimates on $U$ and $V$ obtained by (\ref{eqn:upper_bound_Un}) and 
(\ref{eqn:upper_bound_nabla_Un}). Moreover,
\[ \lim_n \varphi_{t_n}(s) = \varphi_{t}(s) 
   \sptext{1}{for all}{1}
   s\in (r,R)\backslash \{ t \}. \]
As for $\E[ G(y; X_R^{t,x}) N_R^{t,1,(t,x)}]$, we show that 
$\lim_n \psi_{t_n,x_n}(s) = \psi_{t,x}(s)$ for all $s\in (r,R)\setminus \{ t \}$.
\bigskip

{\bf (h)}
Finally, we show that $\nabla_x U = V$. For $x_0,x_1\in \R^d$ we have that
\[   U^n(y;t,x_0) - U^n(y;t,x_1) 
   = \int_0^1 \langle \nabla_x U^n(y;t, x_0+\lambda (x_1-x_0)),x_1-x_0 \rangle d\lambda \]  
for $r\le t < R$.
By dominated convergence we have that
\[   U(y;t,x_0) - U(y;t,x_1) 
   = \int_0^1 \langle V (y;t, x_0+\lambda (x_1-x_0)),x_1-x_0 \rangle d\lambda \]
so that we are done.
\end{proof}
\bigskip

\begin{lemma}
[$L_p$-bound for the $Z$-process for a singular generator]
\label{lemma:upper_bound_V}
Assume condition $(A_{b,\sigma})$, $0\le r < R \le T$, $2\le p < \infty$
and assume that $X=(X_s)_{s\in [r,R]}$ is the diffusion with parameters $(b,\sigma)$ started in some $x_r\in\R^d$.
\footnote{We would need to write $X_s^{r,x_r}$ but use simply $X_s$ to shorten the notation.}
Consider the BSDE
\begin{equation}\label{eqn:linear_BSDE}
U_t = \int_t^R h(s,X_s,U_s,V_s) ds - \int_t^R V_s dB_s
\end{equation}
with a generator 
$h:[r,R)\times \R^d\times \R^d \times \R^{d\times d} \to \R^d$ which is 
measurable with respect to $\mathcal{B}([r,R)) \times \mathcal{B}(\R^d) \times \mathcal{B}(\R^d) \times \mathcal{B}(\R^{d\times d})$
and assume the following:
\begin{enumerate}[{\rm (i)}]
\item $h(s,\cdot,u,v)$ is continuous in $x$ for fixed $s,u,v$.
\item $|h(s,x,u_1,v_1)-h(s,x,u_2,v_2)| \le L (|u_1-u_2| + |v_1-v_2|)$ for some $L>0$.
\item $|h(s,x,u,v)| \le \alpha(s,x) + \lambda |u| + \mu |v|$ where
      $\alpha:[r,R)\times \R^d\to \R$ is non-negative and 
      $\mathcal{B}([r,R)) \times \mathcal{B}(\R^d)$-measurable,
      $\alpha(s,\cdot)$ is continuous for fixed $s$ and
      satisfies $\alpha(s,x) \le \kappa(s)[1+|x|^q]$ for
      some $q\ge 0$, where the  function $\kappa(.)\geq 0$ is bounded on compact subintervals of $[r,R)$ and
      \[ \int_r^R \|\alpha(s,X_s)\|_p ds < \infty.\]
\end{enumerate}
Then there exists an unique solution $(U,V)$ such that
\[ \sup_{r\le t \le R} |U_t| + \kla \int_r^R |V_t|^2 dt \mer^\frac{1}{2} \in L_p \]
and a constant $c=c(p,\sigma,b,T,L,\lambda,\mu)>0$ such that
\begin{enumerate}[{\rm (1)}]
\item $\| U_t \|_p \le c \| \int_t^R |\alpha(s,X_s)| ds \|_p$ for $t\in [r,R)$,
\item and there exists a Borel set $\mathcal{N}\subseteq [r,R)$ of 
      Lebesgue measure zero such that
      \[ \| V_t \|_p \le c \int_t^R \frac{\|\alpha(s,X_s)\|_p}
                                              {\sqrt{s-t}} ds
                                               \]
      for all $t\in [r,R)\backslash\mathcal{N}$.
\end{enumerate}
\end{lemma}

\begin{proof}
The local boundedness of $\kappa$ ensures 
$\int_t^R \frac{\|\alpha(s,X_s)\|_p}  {\sqrt{s-t}} ds<\infty$ for $t\in[r,R)$.
The existence of the unique $L_p$-solution $(U,V)$ follows from  
\cite[Theorems 4.1 and 4.2]{Briand-et-al} and the statement (1) 
follows from \cite[Proposition 3.2]{Briand-et-al} where we consider the BSDE 
with the generator $h^{(t)}(s,x,u,v):= h(s,x,u,v)$ if $s\in [t,R)$ and $h^{(t)}:=0$
otherwise, and the accordingly modified $\alpha$.
So we turn to the statement (2).
\medskip

(a) Fix a bump-function $v:\R^d\to [0,\infty)\in C^\infty_0$ with $v(x)=0$ for $|x|\ge 1$ and 
$\int_{\R^d} v(x)dx=1$. For $N\ge 1$, $\vare>0$, $x\in \R^d$ and $\xi\in\R$ define
\equa
 v_\vare(x) &:= & \frac{1}{\vare^d}v \kla \frac{x}{\vare} \mer, \\ 
 h_{\vare,N} (s,x,u,v) 
             &:= & (v_\vare^x \ast h^N)(s,x,u,v)
\tion
where $h^N:=(h_1^{N/\sqrt{d}},...,h_d^{N/\sqrt{d}})$
with $\xi^N=(\xi\wedge N)\vee (-N)$ for $\xi\in\R$
(so that $|h^N| \le N$) and 
the notation $v_\vare^x$ indicates that the convolution is taken with
respect to $x$. 
Assumption (ii) implies that
\[         |h_{\vare,N}(s,x,u,v)|
      \le  (v_\vare^x \ast \alpha^N)(s,x) + \lambda |u| + \mu |v|. \]
The function $h_{\vare,N}$ is uniformly Lipschitz in $(x,u,v)$ as
\begin{multline*}     
                |h_{\vare,N} (s,x_1,u_1,v_1) - h_{\vare,N} (s,x_2,u_2,v_2) | \\
       \le  L (|u_1-u_2| + |v_1-v_2|) 
                     + \sup_{s',x',u',v'} \bet \nabla_{x'} h_{\vare,N}(s',x',u',v')\rag |x_1-x_2|,
\end{multline*}
where we note that $\nabla_{x'} h_{\vare,N}$ is a matrix,
       and 
       \equa 
           \bet \nabla_{x'} h_{\vare,N}(s',x',u',v')\rag 
       & = & \vare^{-d-1} \bet \int_{|\xi-x'|\le \vare} (\nabla v)\kla \frac{x'-\xi}{\vare}\mer  
             h^N(s,\xi,u',v') d \xi \rag \\
       &\le& \vare^{-1} vol(B_1(\R^d)) N \|\nabla v\|_\infty.
      \tion

(b) Fix $N\ge 1$ and $\vare>0$, let
\[ h_0(s,x,u,v) :=  h_{\vare,N}(s,x,u,v)\chi_{[r,R) }(s) \]
and
$\alpha_0(s,x) := (v_\vare^x \ast \alpha^N)(s,x)\chi_{[r,R) }(s)$. Let $(U^0,V^0)$ be the solution of our BSDE (\ref{eqn:linear_BSDE}) with $h$ replaced by $h_0$
according to \cite[Theorem 2.6]{Delarue:SPA-2002}, where $U^0_s := A^0(s,X_s)$ for a 
continuous and bounded function 
$A^0:[r,R]\times \R^d \to \R^d$. It is also shown that 
$\lambda\times \Q(\{ (t,\omega) \in [r,R]\times M : |V_s^0| > c \}) = 0$ for some 
$c>0$. By considering a Picard iteration 
\[ U^{0,k}_t = \int_t^R h_0(s,X_s,A^0(s,X_s),V^{0,k-1}_s) ds -
               \int_t^R V^{0,k}_s dB_s \]
with $U^{0,0}_s\equiv 0$ one can show by induction that $V^{0,k}_s$ can be realized as a measurable functional 
of $s$ and $X_s$ and obtains finally that   
there is a measurable function $B^0 : [r,R]\times \R^d \to \R^{d\times d}$ with $\| B^0 \|_\infty \le c$,
such that one can realize (using uniqueness from \cite[Theorem 2.6]{Delarue:SPA-2002})
$V^0$ as $V_s^0 = B^0(s,X_s)$. Now
\[ h_0(s,X_s,U_s^0,V_s^0) = \E h_0(s,X_s,U_s^0,V_s^0) + \int_r^s \lambda_t^s dB_t \quad a.s.,\]
where the matrix $\lambda_t^s$ is obtained via the PDE approach, so that we get, a.s.,
\equa
      U_r^0 + \int_r^R V_t^0 dB_t
& = & \int_r^R h_0(s,X_s,U_s^0,V_s^0) ds \\
& = & \int_r^R \E h_0(s,X_s,U_s^0,V_s^0) ds + \int_r^R \int_t^R \lambda_t^s ds dB_t
\tion
by a stochastic Fubini argument
and
\[ V_t^0 = \int_t^R \lambda_t^s ds \quad  a.s.
   \sptext{.3}{for a.e.}{.3}
   t\in [r,R]. \]
If the set of those $t$ is denoted by $\mathcal{M}$, then for $t\in \mathcal{M}$,
\equa
      \| V_t^0 \|_p 
&\le& \int_t^R \|\lambda_t^s \|_p ds \\
&\le& \kappa_{p'} \int_t^R \frac{\|h_0(s,X_s,U_s^0,V_s^0)\|_p}{\sqrt{s-t}} ds\\
&\le& \kappa_{p'} \int_t^R \frac{\|a_0(s,X_s)\|_p+\lambda \|U_s^0\|_p + \mu \|V_s^0\|_p}
                             {\sqrt{s-t}} ds \\
& = & \kappa_{p'} \int_t^R \frac{\psi(s)+ \mu \|V_s^0\|_p}
                             {\sqrt{s-t}} ds                  
\tion
with
\[ \psi(s) := \|\alpha_0(s,X_s)\|_p+\lambda \|U_s^0\|_p.\]
Applying the same inequality to $s\in \mathcal{M}$ gives by iteration for $t\in \mathcal{M}$ that
\equa
&   &  \| V_t^0 \|_p \\
&\le& \kappa_{p'} \int_t^R \frac{\psi(s)+ \mu \kappa_{p'} \int_s^R \frac{\psi(w)+ \mu \|V_w^0\|_p}
                             {\sqrt{w-s}} dw}{\sqrt{s-t}} ds \\
& = &  \kappa_{p'} \int_t^R \frac{\psi(s)}{\sqrt{s-t}} ds
     + \mu \kappa_{p'}^2   B\kla \frac{1}{2},\frac{1}{2}\mer \int_t^R \psi(s) ds + \\
&   & \hspace*{16em}
     + (\mu \kappa_{p'})^2  B\kla \frac{1}{2},\frac{1}{2}\mer \int_t^R \|V_s^0\|_p ds \\
&\le&  \left (\kappa_{p'} + \sqrt{T} \mu \kappa_{p'}^2   B\kla \frac{1}{2},\frac{1}{2}\mer \right ) 
        \int_t^R \frac{\psi(s)}{\sqrt{s-t}} ds \\
&   & \hspace*{15em} + (\mu \kappa_{p'})^2  B\kla \frac{1}{2},\frac{1}{2}\mer \int_t^R \|V_s^0\|_p ds.
\tion
It follows from the boundedness properties of $V_s^0$ that 
$\int_r^R \|V_s^0\|_p ds < \infty$. For this reason we can 
apply Gronwall's lemma to derive
\[
    \| V_t^0 \|_p 
\le (\kappa_{p'} + \sqrt{T} \mu \kappa_{p'}^2   B(1/2,1/2))  
    e^{(\mu \kappa_{p'})^2 B(1/2,1/2)(R-t)}
    \int_t^R \frac{\psi(s)}{\sqrt{s-t}} ds
\]
for $t\in \mathcal{M}$.
Next we estimate $\psi(s)$ by
\[      \psi(s) 
   \le \| \alpha_0(s,X_s) \|_p + \lambda c_{(\ref{lemma:upper_bound_V})(1)}
       \int_s^R \|\alpha_0(w,X_w) \|_p dw, \]
where we use Lemma \ref{lemma:upper_bound_V}(1) (with the same 
$(L,\lambda,\mu)$), and get
\equa
      \| V_t^0 \|_p 
&\le& d_1 \int_t^R \frac{\|\alpha_0(s,X_s)\|_p+\int_s^R \|\alpha_0(w,X_w) \|_p dw}{\sqrt{s-t}} ds\\
&\le&  d_1 (1+2T) \int_t^R \frac{\|\alpha_0(s,X_s)\|_p}{\sqrt{s-t}} ds
\tion
with 
$d_1:= \left (\kappa_{p'} + \sqrt{T} \mu \kappa_{p'}^2   B\kla \frac{1}{2},\frac{1}{2}\mer \right )  
    e^{(\mu \kappa_{p'})^2  B\kla \frac{1}{2},\frac{1}{2}\mer  T} (1+\lambda c_{(\ref{lemma:upper_bound_V})(1)})$.
Hence, re\-writing the dependence with respect to $N$ and $\epsilon$ in our estimates, we have proved
\begin{equation}\label{eqn:upper_bound_V}
\| V_t^{N,\vare} \|_p \le d_2 \int_t^R \frac{\|(v_\vare^x \ast \alpha^N)(s,X_s)\|_p}
                                              {\sqrt{s-t}} ds
\end{equation}
for $t \in \mathcal{M}=[r,R]\backslash \mathcal{N}_{N,\vare}$ with $d_2:= d_1 (1+2T)$.
\smallskip

(c) Let $\overline{\mathcal{N}}:= \bigcup_{N,n} \mathcal{N}_{N,1/n}$ and let
$((U_t^N,V_t^N))_{t\in [r,R]}$ be the solution of (\ref{eqn:linear_BSDE}) with the generator $h^N$. Because 
\[ \lim_{\vare\downarrow 0} \int_r^R 
    \| h_{N,\vare}(s,X_s,U_s^N,V_s^N) -h^N(s,X_s,U_s^N,V_s^N) \|_2 ds = 0 \]
by dominated convergence (here we use the continuity of $h$ in $x$) and
\[ |h_{N,\vare}(r,x,u_1,v_1) - h_{N,\vare}(r,x,u_2,v_2) |\le L [|u_1-u_2| + |v_1-v_2|], \]
Lemma \ref{lemma:compare_briand_etal} implies that
\[  \lim_{n\to \infty} \int_r^R 
    \| V_s^{N,1/n} -V^N_s\|_2^2 ds = 0 \]
for all $N=1,2,...$
Hence there are sub-sequences $(n_l^N)_{l=1}^\infty$ such that
\[ \lim_{l\to \infty} 
    | V_s^{N,1/{n_l^N}} -V^N_s| = 0
    \quad \Q \times \lambda \quad a.e.\]
and a Borel set $\mathcal{N}_N\subseteq [0,T]$ of
Lebesgue measure zero such that 
\[ V_s^{N,1/{n_l^N}} \to_l V_s^N \sptext{1}{a.s. for}{1} s\not\in \mathcal{N}_N. \] 
Applying Fatou's lemma on the left-hand side of (\ref{eqn:upper_bound_V})
and dominated convergence on the right-hand side (note that
$|v^x_\vare\ast \alpha^N|\le N$ and that $\alpha$ is supposed to be continuous in $x$), we derive
\[ \| V_t^{N} \|_p \le d_2 \int_t^R \frac{\|\alpha^N(s,X_s)\|_p}
                                              {\sqrt{s-t}} ds
                   \le d_2 \int_t^R \frac{\|\alpha(s,X_s)\|_p}
                                              {\sqrt{s-t}} ds
                                               \]
for all $t\in [r,R]\backslash (\bigcup_{N'=1}^\infty \mathcal{N}_{N'} \cup\overline{\mathcal{N}})$.
In the same way, Lemma \ref{lemma:compare_briand_etal},
\[ \int_r^R \| h^N(s,X_s,U_s,V_s) -h(s,X_s,U_s,V_s) \|_2 ds \to_N 0 \]
and 
$|h^N(s,x,u_1,v_1) - h^N(s,x,u_2,v_2) |\le L [|u_1-u_2| + |v_1-v_2|]$
give 
\[  \lim_{N\to \infty} \int_r^R 
    \| V_s^{N} -V_s\|_2^2 ds = 0 \]
and the existence of a subsequence $(N_k)_{k=1}^\infty$ such that
\[ \lim_{k\to \infty} 
    | V_s^{N_k} -V_s| = 0
    \quad \Q \times \lambda \quad a.e.\]
Hence there is some $\mathcal{N}_0\subseteq [r,R]$ of Lebesgue measure
zero such that
\[ V_s^{N_k} \to_k V_s \sptext{1}{a.s. for}{1} s\not\in \mathcal{N}_0. \] 
Again applying Fatou's lemma gives that
\[ \| V_t \|_p \le d_2 \int_t^R \frac{\|\alpha(s,X_s)\|_p}
                                              {\sqrt{s-t}} ds
                                               \]
for all $N=1,2,...$ and  $t\in [r,R]\backslash (\bigcup_{N'=0}^\infty \mathcal{N}_{N'}\cup\overline{\mathcal{N}})$.
\end{proof}

\section{Appendix}

\begin{proposition}
[{\cite[pp. 260, 72, 74, 44] {Friedman-0}}]
\label{proposition:transition_density}
For $b,\sigma$ satisfying ($A_{b,\sigma}$), there exists a continuous transition density 
\[ \Gamma: \{(t,x,s,\xi):0\leq t<s\leq T \mbox{ and } x,\xi \in\R^d\}\to (0,\infty) \]
such that
$\P(X^{t,x}_s\in B) = \int_B \Gamma(t,x;s,\xi)d\xi$ for $0\leq t<s \leq T$ and
$B\in \mathcal{B}(\R^d)$, where
\[ X^{t,x}_s = x + \int_t^s b(r,X^{t,x}_r)dr + \int_t^s \sigma(r,X^{t,x}_r) dW_r, \]
such that the following is satisfied:
\begin{enumerate}[{\rm (i)}]
\item For all multi-indices $m$ and $k$ with $|m|+2k \le 3$ one has that the derivatives
      $D_t^k D_x^m \Gamma(t,x;s,\xi)$ exists and are continuous on $[0,s)\times \R^d$ 
      and that the differentiation can be done in any order.
\item For $0\leq t<s \leq T$ and $(x,\xi)\in\R^d\times\R^d$ one has 
      \[    \frac{\partial}{\partial t} \Gamma 
                + \frac{1}{2}  \langle A , D^2 \Gamma \rangle 
                + \langle b,\nabla_x \Gamma \rangle = 0 \]
     where $A= \sigma \sigma^*$ and 
     $D^2 =  \left ( \frac{\partial^2}{\partial x_i \partial x_j} \right )_{i,j=1}^d$.
\item For all multi-indices $m$ with $|m|\le 3$ there exists a constant $c=c_m>0$ such that for
      $0\leq t<s \leq T$ and $(x,\xi)\in\R^d\times\R^d$ one has that
      \[      \bet D^m_x \Gamma (t,x;s,\xi)\rag 
          \le c \pl (s-t)^{-\frac{|m|}{2}} \gamma^{d}_{s-t}\kla\frac{x-\xi}{c}\mer \]
          
      where 
      $   \gamma_t^{d}(\eta) 
       := \frac{1}{(2\pi  t)^{d/2}}
           e^{-\frac{|\eta|^2}{2t}}$.
\end{enumerate}
\end{proposition}

\begin{remark}\label{rem:Malliavin_weights}
\rm
The weights $N^{r,i,(t,x)}_R$ are essential so that we briefly recall their construction.
For notational simplicity we let $t=0$ and omit the superscripts $(t,x)$. For $i=1$ one has 
$N_R^{r,1}:= \frac{1}{R-r} \kla \int_r^R  (\sigma(s,X_s)^{-1} \nabla X_s\nabla X_r^{-1})^* dW_s \mer^*$
where $\nabla X_t = \nabla_x b(t,X_t) \nabla X_t dt +  \nabla_x \sigma(t,X_t) \nabla X_t dW_t$ with
$\nabla X_0=I_{\R^d}$, the identity matrix (see, for example, \cite{Ma-Zhang:2002,gobe:munos:05}). 
To consider $i=2$ we follow \cite{gobe:munos:05} and 
let $0\le r < R \le T$, $\rho:= (r+R)/2$, $g:\R^d\to \R$ be a Borel measurable
polynomially bounded function and $F$ like in (\ref{eqn:F}). For $k=1,...,d$ we have
that
$(\partial F/\partial x_k)(r,X_r) = \E (F(\rho,X_\rho) N_\rho^{r,1}(k)|\ftn_r)$ a.s.
Applying the $\nabla$-operator, which can be justified by standard methods, we derive that, a.s.
\equa
&   &   \nabla_x (\partial F/\partial x_k)(r,X_r) \nabla X_r \\
& = &   \E (\nabla_x F(\rho,X_\rho) \nabla X_\rho N_\rho^{r,1}(k)
      + F(\rho,X_\rho) \nabla N_\rho^{r,1}(k) |\ftn_r) \\
& = &   \E ( \E (g(X_R)N^{\rho,1}_R|\ftn_\rho) \nabla X_\rho N_\rho^{r,1}(k)
      + \E (g(X_R)|\ftn_\rho)  \nabla N_\rho^{r,1}(k) |\ftn_r).
\tion
Therefore we can take
$N_R^{r,2}(k) :=  [   N^{\rho,1}_R \nabla X_\rho N_\rho^{r,1}(k)
                    + \nabla N_\rho^{r,1}(k)] (\nabla X_r)^{-1}$
to obtain
$ \nabla_x (\partial F/\partial x_k)(r,X_r) = \E (g(X_R)N_R^{r,2}(k)|\ftn_r)$ a.s.
\end{remark}

\paragraph{Proof of Proposition \ref{proposition:change_gamma}.}
Assume that we have the diffusions $X^1 =(X_t^1)_{t\in [0,T_1]}$ and 
$X^2 =(X_t^2)_{t\in [0,T_2]}$ starting in $x_1\in \R^d$ and $x_2\in \R^d$ respectively,
satisfying our assumptions with the corresponding 
transition densities $\Gamma_1$ and $\Gamma_2$, and assume that they satisfy 
\[ \Gamma_1(t,x;s,\xi) \le M \Gamma_2(\mu t, \nu x; \mu s,\nu \xi) \]
for some $M,\mu,\nu>0$ and all $x,\xi \in \R^d$ and $0\le t  \le s \le T_1$
and with $T_2=\mu T_1$.
Let $g:\R^d\to \R$ be a polynomially bounded Borel function.
Then, for $x_2=\nu x_1$, 
\[      \E| g(X_{T_1}^1) - \E (g(X_{T_1}^1)|\ftn_t)|^p
    \le 2^p \frac{M^3}{\nu^3} \E| \widetilde{g}(X_{\mu T_1}^2) - \E (\widetilde{g}(X_{\mu T_1}^2)|\ftn_{\mu t})|^p
\]
with
\[ \widetilde{g}(x):= g\left (\frac{x}{\nu}\right ).\]
In fact, we have that
\equa
&   & \E| g(X_{T_1}^1) - \E (g(X_{T_1}^1)|\ftn_t)|^p \\
&\le& \int_{\R^d} \int_{\R^d} \int_{\R^d}
      |g(\xi) - g(\eta)|^p \\
&   & \hspace*{5em}      \Gamma_1(0,x_1;t,x) \Gamma_1(t,x;T_1,\xi) \Gamma_1(t,x;T_1,\eta) d x d\xi d\eta \\
&\le& M^3 \int_{\R^d} \int_{\R^d} \int_{\R^d}
      |g(\xi) - g(\eta)|^p \\
&   & \hspace*{5em}      \Gamma_2(0,\nu x_1;\mu t,\nu x) \Gamma_2(\mu t,\nu x;\mu T_1,\nu \xi) 
                         \Gamma_2(\mu t,\nu x;\mu T_1,\nu \eta) d x d\xi d\eta \\
& = & \frac{M^3}{\nu^3} \int_{\R^d} \int_{\R^d} \int_{\R^d}
      |\widetilde{g}(\xi) - \widetilde{g}(\eta)|^p \\
&   & \hspace*{5em}      \Gamma_2(0, x_2;\mu t,x) \Gamma_2(\mu t,x;T_2,\xi) \Gamma_2(\mu t,x;T_2,\eta) d x d\xi d\eta \\
&\le& 2^{p} \frac{M^3}{\nu^3} 
      \E| \widetilde{g}(X_{T_2}^2) - \E (\widetilde{g}(X_{T_2}^2)|\ftn_{\mu t})|^p
\tion
where we used relation (\ref{eqn:conditional_expectation_vs_transition_density}).
This implies our assertion by taking
$(\Gamma_1,x_1,T_1)=(\Gamma,x_0,r_l)$ and $T_2=T_1$, $\nu=1/c_{(\ref{proposition:transition_density})}$ and
$X_t^2=\nu x_0+W_t$.
\hfill$\Box$


\end{document}